%% file: SpaceTimeINS.tex
\documentclass[3p,times]{elsarticle}

\usepackage{bbm}
\usepackage{bm}
\usepackage{color}
\usepackage{xcolor}

\usepackage{graphicx}
\usepackage{epstopdf}
\usepackage{psfrag}

\usepackage{amssymb}
\usepackage{amsmath}
\usepackage{dsfont}
\usepackage{rotating}
\usepackage{pgf,tikz}
\usetikzlibrary{arrows}

\definecolor{ttzzqq}{rgb}{0.2,0.6,0}

\definecolor{qqttcc}{rgb}{0,0.2,0.8}

\definecolor{qqttzz}{rgb}{0,0.2,0.6}
\definecolor{ffqqqq}{rgb}{1,0,0}
\definecolor{qqwuqq}{rgb}{0,0.39,0}
\definecolor{zzttqq}{rgb}{0.6,0.2,0}
\definecolor{qqqqff}{rgb}{0,0,1}
\definecolor{ttttqq}{rgb}{0.2,0.2,0}

\definecolor{qqwwtt}{rgb}{0,0.4,0.2}
\definecolor{ubqqys}{rgb}{0.29,0,0.51}

\definecolor{wwttqq}{rgb}{0.4,0.2,0}
\definecolor{uuuuuu}{rgb}{0.27,0.27,0.27}
\definecolor{qqzzff}{rgb}{0,0.6,1}
\definecolor{xdxdff}{rgb}{0.49,0.49,1}
\definecolor{ccwwqq}{rgb}{0.8,0.4,0}

\definecolor{ttqqqq}{rgb}{0.2,0,0}

\definecolor{qqzzcc}{rgb}{0,0.6,0.8}

\newcommand{\R}{\mathbb{R}}

\newcommand{\eref}[1]{$(\ref{#1})$}
\newcommand{\p}{\mathbb{\wp}}
\newcommand{\np}{p}

\newcommand{\nextud}{\vec{n}_{ij}}
\newcommand{\nstd}{\vec{n}_{j}}
\newcommand{\etah}{ \hat{\bm{p}}}
\newcommand{\bphi}{ \bm{\phi}}
\newcommand{\bpsi}{ \bm{\psi}}

\newcommand{\vvh}{\hat{\bm{\mathbf{v}}}}

\newcommand{\Fvh}{\widehat{\bm{F\mathbf{v}}}}

\newcommand{\nv}{\vec{n}}
\newcommand{\B}{\mathcal{B}}
\newcommand{\TF}{F}

\newcommand{\TT}{\mbox{\boldmath$T$}}
\newcommand{\QQ}{\mbox{\boldmath$R$}}
\newcommand{\poh}{+\frac{1}{2}}

\newcommand{\st}{ {st}}
\newcommand{\Ss}{ \bm{\mathcal{S}}}

\newcommand{\D}{\bm{\mathcal{D}}}
\newcommand{\Q}{\bm{\mathcal{Q}}}
\newcommand{\RM}{\bm{\mathcal{R}}}
\newcommand{\LM}{\bm{\mathcal{L}}}

\newcommand{\Mpsi}{\bm{M}}

\newcommand{\diff}[2]{\frac{\partial {#1} }{\partial {#2} } }

\usepackage{amssymb}

\journal{Journal of computational physics}

\begin{document}

\begin{frontmatter}



\title{A staggered space-time discontinuous Galerkin method for the incompressible Navier-Stokes equations on two-dimensional triangular meshes} 

\author[1]{Maurizio Tavelli\fnref{label1}}
\author[2]{Michael Dumbser \corref{corr1} \fnref{label2}}
\address[1]{Department of Mathematics, University of Trento, Via Sommarive 14, I-38050 Trento, Italy}
\address[2]{Department of Civil, Environmental and Mechanical Engineering, University of Trento, Via Mesiano 77, I-38123 Trento, Italy}

\fntext[label1]{\tt m.tavelli@unitn.it (M.~Tavelli)}
\fntext[label2]{\tt michael.dumbser@unitn.it (M.~Dumbser)}
\begin{abstract}
In this paper we propose a novel arbitrary high order accurate semi-implicit space-time discontinuous Galerkin method for the solution of the two dimensional incompressible Navier-Stokes equations on \textit{staggered}  unstructured triangular meshes. Isoparametric finite elements are used to take into account curved domain boundaries. The discrete pressure is defined on the primal triangular grid and the discrete  velocity field is defined on an edge-based staggered dual grid. While staggered meshes are state of the art in classical finite difference approximations of the incompressible Navier-Stokes equations, their use in the context of high order DG schemes is novel and still quite rare. 
Formal substitution of the discrete momentum equation into the discrete continuity equation yields a sparse four-point block system for the scalar pressure, which is conveniently solved with a matrix-free
GMRES algorithm. A very simple and efficient Picard iteration is then used in order to achieve high order of accuracy also in time, which is in general a non-trivial task in the 
context of high order discretizations for the incompressible Navier-Stokes equations. 
The flexibility and accuracy of high order space-time DG methods on curved unstructured meshes allows to discretize even complex physical domains with very coarse grids in both, space and time. The 
use of a staggered grid allows to \textit{avoid} the use of Riemann solvers in several terms of the discrete equations and significantly reduces the total stencil size of the linear system  
that needs to be solved for the pressure. 
The proposed method is validated for approximation polynomials of degree up to $p=4$ in space and time by solving a series of typical numerical test problems and by comparing the obtained numerical 
results with available exact analytical solutions or other numerical reference data.  
\end{abstract}


\end{frontmatter}


\section{Introduction}
The discretization of the incompressible Navier-Stokes equations was mainly carried out in the past using finite difference methods \cite{markerandcell,patankarspalding,patankar,vanKan} as well as  continuous finite element schemes \cite{TaylorHood,SUPG,SUPG2,Fortin,Verfuerth,Rannacher1,Rannacher3}. 
On the contrary, the construction of high order discontinuous Galerkin (DG) finite element methods for the incompressible Navier-Stokes equations is still a very active topic of ongoing research. 
Obtaining high order of accuracy also in time represents an important goal in order to achieve accurate results for unsteady problems. 

Several high order DG methods for the incompressible Navier-Stokes equations have been recently presented in the literature, see for example     
\cite{Bassi2007,Shahbazi2007,Ferrer2011,Nguyen2011,Rhebergen2012,Rhebergen2013,Crivellini2013,KleinKummerOberlack2013}, without pretending completeness. 
An alternative is the DG scheme proposed by Bassi et al. in \cite{Bassi2006}, which is based on an extension of the technique of artificial compressibility that was 
originally introduced by Chorin in the finite difference context \cite{chorin1,chorin2}. 
Another very well  known approach to discretize general convection-diffusion equations in the context of hp discontinuous Galerkin finite element methods is the one  
proposed by Baumann and Oden in \cite{Baumann1999311,Baumann199979}. A unified analysis of several variants of the DG method applied to an 
elliptic model problem has been provided by Arnold et al. in \cite{ArnoldBrezzi}. We also would like to mention recent works on semi-implicit DG schemes, 
such as the ones presented in \cite{TumoloBonaventuraRestelli,GiraldoRestelli,Dolejsi1,Dolejsi2,Dolejsi3}, to which our approach is indirectly related. 

While the use of staggered grids is a very common practice in the finite difference community, its use is not so widespread in the development of high order DG schemes. 
The first staggered DG schemes, based on a \textit{vertex-based} dual grid, have been proposed in  \cite{CentralDG1,CentralDG2}. Other recent high order staggered DG schemes that 
use an \textit{edge-based} dual grid have been forwarded in \cite{StaggeredDGCE1,StaggeredDG,DumbserCasulli}. 
The advantage in using edge-based staggered grids is that they allow to improve significantly the sparsity pattern of the final linear system that has to be solved for 
the pressure.  

Very recently, a staggered semi-implicit DG scheme for the solution of the two dimensional shallow water equations was presented in \cite{DumbserCasulli,2DSIUSW} and then extended in \cite{2STINS} to 
the incompressible Navier-Stokes equations. The method presented in \cite{2STINS} is in principle of arbitrary high order of accuracy in space, while it reaches only second order 
in time. Consequently, it does not allow to recover high order accurate results for fully unsteady solutions. 


In this paper we propose a new method that is based on the general ideas put forward in \cite{DumbserCasulli,2DSIUSW,2STINS}, but which is also able to reach 
high order of accuracy in time. For this purpose we construct an arbitrary high order accurate \textit{staggered space-time} discontinuous Galerkin finite 
element scheme. By relying on staggered grids we follow the classical philosophy of staggered finite difference schemes for the incompressible Navier-Stokes 
equations and for the free surface shallow water and Navier-Stokes equations, see 
\cite{markerandcell,patankarspalding,patankar,vanKan,HirtNichols,CasulliCheng1992,Casulli1999,CasulliWalters2000,Casulli2009,CasulliStelling2011,CasulliVOF}. 
In the context of staggered finite difference schemes we also would like to mention the so-called multiple pressure variables approach (MPV) 
\cite{klein,RoMu99,ParkMPV}, which is based on the asymptotic analysis of the compressible Navier-Stokes equations and is able to preserve also their incompressible 
limit. 

Our staggered semi-implicit space-time DG method proposed in this paper can be seen as a natural extension of the staggered semi-implicit DG scheme proposed in \cite{2STINS} to 
arbitrary high order of accuracy also in time. However, we emphasize that this extension is not straightforward for the complete convective-viscous problem. 
In the staggered DG scheme presented in \cite{2STINS}, the discrete pressure is defined on the control volumes of the primal triangular mesh, while the discrete 
velocity vector field is defined on an edge-based, quadrilateral dual mesh. In the proposed staggered space-time DG scheme, the spatial control volumes are simply 
extended to the corresponding space-time control volumes by using the tensor product of the spatial control volume with the time interval of each time step, hence 
leading to triangular base prisms for the primal mesh and to quadrilateral base prisms for the dual mesh. 

The nonlinear convective terms are discretized explicitly by using a standard DG scheme \cite{cbs4,CBS-convection-diffusion,CBS-convection-dominated} based on the local 
Lax-Friedrichs (Rusanov) flux \cite{Rusanov:1961a}, while the viscous terms are discretized implicitly using a fractional step method \footnote{Note that high order in time
is obtained later by a Picard iteration, which embraces the entire scheme in each time step.}. The DG discretization of the 
viscous fluxes is based on the formulation of Gassner et al. \cite{MunzDiffusionFlux}, who obtained the viscous numerical flux from the solution of the Generalized 
Riemann Problem (GRP) of the diffusion equation. 

The discrete momentum equation is then inserted into the discrete continuity equation in order to
obtain the discrete form of the pressure Poisson equation. The chosen dual grid used here is taken as the one used in
\cite{Bermudez1998,Bermudez2014,USFORCE,StaggeredDG,2DSIUSW}, which leads to a sparse four-point block system for the scalar pressure. 
Once the new pressure field is known, the velocity vector field can subsequently be updated directly.
A Picard iteration procedure that embraces the entire scheme in each time step is then used in order to achieve arbitrary high-order of accuracy also in time for the 
nonlinear convective term, without introducing a non linearity in the system for the pressure. 


The rest of the paper is organized as follows: in Section \ref{sec_1} the numerical method is described in detail, while in Section \ref{sec.tests}
a set of numerical test problems is solved in order to study the spatial and temporal accuracy of the presented approach. Some concluding remarks are given in Section
\ref{sec.concl}.

\section{DG scheme for the 2D incompressible Navier-Stokes equations}
\label{sec_1}
\subsection{Governing equations}
The two dimensional incompressible Navier-Stokes equations are given by
\begin{eqnarray}
    \frac{\partial \mathbf{v}}{\partial t}+\nabla \cdot \mathbf{F}_c + \nabla p & = & \nu \Delta \mathbf{v} + \mathbf{S} \label{eq:CS_2_2_0}, \\ 
    \nabla \cdot \mathbf{v} & = & 0 \label{eq:CS_2},
\end{eqnarray}
where $p=P/\rho$ indicates the normalized fluid pressure; $P$ is the physical pressure and $\rho$ is the constant fluid density; $\nu$ is the kinematic viscosity  coefficient; $\mathbf{v}=(u,v)$ is the velocity vector; $u$ and $v$ are the velocity components in the $x$ and $y$ direction, respectively;
$\mathbf{S}=\mathbf{S}(\mathbf{v},x,y,t)$ is a (nonlinear) algebraic source term; 
$\mathbf{F}_c=\mathbf{v} \otimes \mathbf{v}$ is the flux tensor of the nonlinear convective terms, namely:
$$ \mathbf{F}_c=\left(\begin{array}{cc} uu & uv \\ vu & vv \end{array} \right). $$

Following the same idea of \cite{MunzDiffusionFlux,ADERNSE}, the viscosity term is first written as $\nu \Delta \mathbf{v}=\nabla \cdot (\nu \nabla \mathbf{v})$ and then grouped with the nonlinear convective term.
In this way the momentum Eq. \eref{eq:CS_2_2_0} can be rewritten as
\begin{equation}
	\frac{\partial \mathbf{v}}{\partial t}+\nabla \cdot \mathbf{\TF} + \nabla p= \mathbf{S} 
\label{eq:CS_2_2},
\end{equation}
where $\mathbf{\TF}=\mathbf{\TF}(\mathbf{v},\nabla \mathbf{v})=\mathbf{F}_c(\mathbf{v})-\nu \nabla \mathbf{v}$ is a nonlinear tensor that depends on the velocity and its gradient,
see e.g. \cite{MunzDiffusionFlux,ADERNSE}.

\subsection{Staggered unstructured grid}
Through this paper we use the same unstructured staggered grid in space as the one used in \cite{2STINS,2DSIUSW}. In the following, we briefly summarize the grid construction and 
the main notation for the spatial grid. After that, the primary and dual spatial elements are extended to the primary and dual space-time control volumes, respectively.

\subsubsection{Unstructured staggered grid in space}
The spatial computational domain is covered with a set of $N_i$ non-overlapping triangles $\TT_i$ with $i=1 \ldots N_i$. By denoting with $N_j$ the total number of edges, the physical $j-$th edge will be called $\Gamma_j$. $\B(\Omega)$ denotes the set of indices $j$ corresponding to boundary edges.
The three edges of each triangle $\TT_i$ constitute the set $S_i$ defined by $S_i=\{j \in [1,N_j] \,\, | \,\, \Gamma_j \mbox{ is an edge of }\TT_i \}$. For every $j\in [1\ldots N_j]-\B(\Omega)$ there exist two triangles $i_1$ and $i_2$ that share $\Gamma_j$. We assign arbitrarily a left and a right triangle called $\ell(j)$ and $r(j)$, respectively. The standard positive direction is assumed to be from left to right. Let $\nv_{j}$ denote the unit normal
vector defined on the edge $j$ and oriented with respect to the positive direction from left to right. For every triangular element $i$ and edge $j \in S_i$,
the neighbor triangle of element $\TT_i$ that share the edge $\Gamma_j$ is denoted by $\p(i,j)$.
\par
\begin{figure}[ht]
    \begin{center}
    \input{./ugrid.tex}
    \caption{Example of a triangular mesh element with its three neighbors and the associated staggered edge-based dual control volumes, together with the notation
    used throughout the paper.}
    \label{fig.1}
		\end{center}
\end{figure}
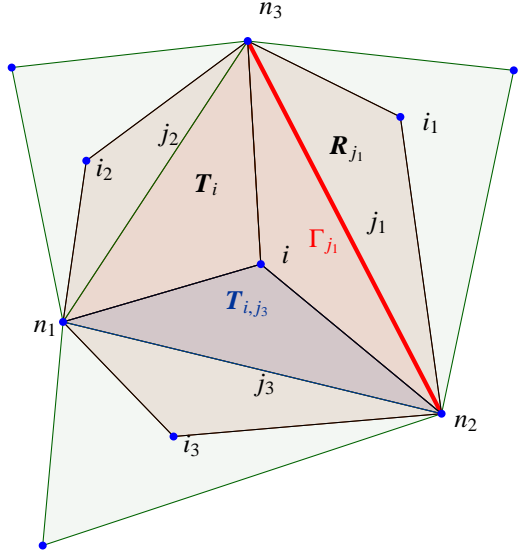

 For every $j\in [1, N_j]-\B(\Omega)$ the quadrilateral element associated to $j$ is called $\QQ_j$ and it is defined, in general, by the two centers of gravity of $\ell(j)$ and $r(j)$ and the two terminal nodes of $\Gamma_j$, see also \cite{Bermudez1998,USFORCE, 2DSIUSW}. We denote by $\TT_{i,j}=\QQ_j \cap \TT_i$ the intersection element for every $i$ and $j \in S_i$. Figure $\ref{fig.1}$ summarizes the used notation, the primal triangular mesh and the dual quadrilateral grid. 
According to \cite{2STINS}, we will call the mesh of triangular elements $\{\TT_i \}_{i \in [1, N_i]}$ the \textit{main grid} or \textit{primal grid} and the quadrilateral grid $\{\QQ_j \}_{j \in [1, N_j]}$ is termed the \textit{dual grid}. 

The dual grid is covering $\Omega$ with non-overlapping quadrilaterals, so we define the equivalent quantities given for the main grid also to the dual one, briefly: $N_l$ is the total amount 
of edges of $\QQ_j$; $\Gamma_l$ indicates the physical $l$-th edge; $\forall j$, the set of edges $l$ of $j$ is indicated with $S_j$; $\forall l$, $\ell_{jl}(l)$ and $r_{jl}(l)$ are the left and the right quadrilateral element, respectively; $\nv_{l}$ is the standard normal vector defined on $l$ and assumed positive with respect to the standard orientation on $l$ (defined, as for the main grid, from the left to the right).
Finally, each triangle $\TT_i$ is defined starting from an arbitrary node and oriented in counter-clockwise direction. Similarly, each quadrilateral element $\QQ_j$ is defined
starting from the point $\ell(j)$ and oriented in counter-clockwise direction.

\subsubsection{Space-time extension}
In the time direction we cover the time interval $[0,T]$ with a sequence of times $0=t^0<t^1<t^2 \ldots <t^N<t^{N+1}=T$. We denote the time step by $\Delta t^{n+1} = t^{n+1}-t^{n} $ and 
the corresponding time interval by $T^{n+1}=[t^{n}, t^{n+1}]$ for $n=0 \ldots N$. In order to ease notation, sometimes we will use the abbreviation $\Delta t= \Delta t^{n+1}$. 
In this way the generic space-time element defined in the time interval $[t^n, t^{n+1}]$ is given by $\TT_i^\st = \TT_i \times T^{n+1}$ for the main grid and $\QQ_j^\st=\QQ_j \times T^{n+1}$ for the dual grid.
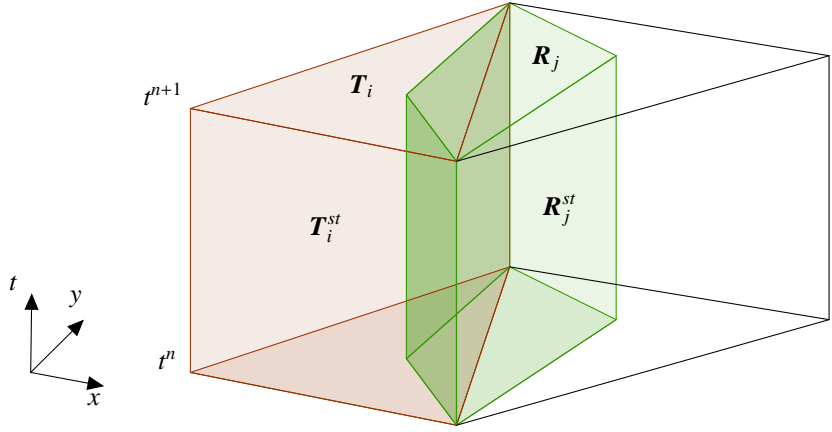
\begin{figure}[ht]
    \begin{center}
	\input{./Omega_st.tex}
    \caption{Example of space-time elements $\TT_i^\st$ (red) and $\QQ_j^\st$ (green) with $j \in S_i$}
    \label{fig.st1}
		\end{center}
\end{figure}
Figure \ref{fig.st1} shows a graphical representation of the primary and dual space-time control volumes. 

\subsection{Space-time basis functions}
\label{sec222}
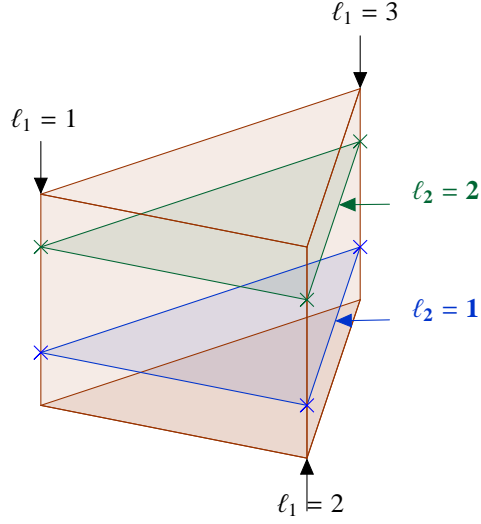
\begin{figure}[ht]
    \begin{center}
		\input{./Layer.tex}
    \caption{Values of the sorting funcitons $\ell_1$ and $\ell_2$ for a triangular element $\TT_i^\st$ with $p=p_\gamma=2$. The cross points represent the DoF in the space-time element}
    \label{fig.st2}
		\end{center}
\end{figure}

According to \cite{2DSIUSW,2STINS} we proceed as follows: we first construct the polynomial basis up to a generic polynomial degree $p$ on some triangular and quadrilateral reference elements. In order to do this we take $T_{std}=\{(\xi,\gamma) \in \R^{2,+} \,\, | \,\, \gamma\leq1-\xi \vee 0 \leq \xi \leq 1 \}$ as the reference triangle and the unit square as the reference quadrilateral element $R_{std}=[0,1]^2$. Using the standard nodal approach of conforming continuous finite elements,  we obtain $N_\phi=\frac{(p+1)(p+2)}{2}$ basis functions $\{\phi_k \}_{k \in [1,N_\phi]}$ on $T_{std}$ and $N_{\psi}=(p+1)^2$ basis functions on $R_{std}$.
The connection between reference and physical space is performed by the maps $T_i:\TT_i \longrightarrow T_{std}$ for every $i =1 \ldots N_i$; $T_j:\QQ_j \longrightarrow R_{std}$ for every $j =1 \ldots N_j$ and its inverse, called $T_i^{-1}:\TT_i \longleftarrow T_{std}$ and $T_j^{-1}:\QQ_j \longleftarrow R_{std}$, respectively. The maps from the physical coordinates to the reference one can be constructed following a classical
sub-parametric or a complete iso-parametric approach.
In the same way we construct the time basis functions on a reference interval $[0,1]$ for polynomials of degree $p_\gamma$. 
In this case the resulting $N_\gamma=p_\gamma+1$ basis functions $\{\gamma_k\}_{k \in [1, N_\gamma]}$ are defined as the Lagrange interpolation polynomials passing through the Gauss-Legendre quadrature points for the unit interval. For every time interval $[t^n, t^{n+1}]$, the map between the reference interval and the physical one is simply given by $t=t^n+\tau \Delta t^{n+1}$ for every $\tau \in [0,1]$.
Using the tensor product we can finally construct the basis functions on the space-time elements $\TT_i^\st$ and $\QQ_j^\st$ such as $\tilde{\phi}(\xi,\gamma,\tau)=\phi(\xi, \gamma) \cdot \gamma(\tau)$ and $\tilde{\psi}(\xi, \gamma,\tau)=\psi(\xi, \gamma) \cdot \gamma(\tau)$. The total number of basis functions becomes $N_\phi^\st=N_\phi \cdot N_\gamma$ and $N_\psi^\st=N_\psi \cdot N_\gamma$. By introducing two sorting functions $\ell_1(\,\, ,N_\cdot^\st):[1,N_\cdot^\st] \rightarrow [1,N_\cdot]$ and $\ell_2(\,\, ,N_\cdot^\st):[1,N_\cdot^\st] \rightarrow [1,N_\gamma]$, defined as
\begin{eqnarray}
    \ell_2(k,N) &=& int\left[\frac{k-1}{N}\right]+1 \nonumber \\
    \ell_1(k,N) &=& k-(\ell_2(k,N)-1)\cdot N
\label{eq:DST0}
\end{eqnarray}

 we can explicit the form of $\tilde{\phi}_k$ and $\tilde{\psi}_l$ for $k=1 \ldots N_\phi^\st$ and $l=1 \ldots N_\psi^\st$ in terms of space and time basis functions:
\begin{eqnarray}
    \tilde{\phi}_k(\xi,\gamma,\tau) &=& \phi_{\ell_1(k,N_\phi^\st)}(\xi,\gamma) \gamma_{\ell_2(k,N_\phi^\st)}(\tau) \qquad \forall k\in [1, N_\phi^\st] \nonumber \\
    \tilde{\psi}_k(\xi,\gamma,\tau) &=& \psi_{\ell_1(k,N_\psi^\st)}(\xi,\gamma) \gamma_{\ell_2(k,N_\psi^\st)}(\tau) \qquad \forall k\in [1, N_\psi^\st] \nonumber
\label{eq:DST1}
\end{eqnarray}
Remark how $\ell_2$ can be seen as a temporal layer selector function, so all the indexes $k$ such that $l_2(k,\cdot)=l$ represent the spatial degrees of freedom (DoF) at the time layer $l$, for every  fixed $l=1 \ldots N_\gamma$. In the same way $l_1(k,\cdot)=m$ represents the time evolution of the DoF $m$ inside the space-time element $\TT_i^\st$. An example of how the sorting functions act is shown in Figure \ref{fig.st2}. 

\subsection{Semi-implicit space-time DG scheme} 
\label{sec_semi_imp_dg}
The discrete pressure $p_h$ is defined on the main grid, namely $p_h(x,y,t)|_{\TT_i^\st}=p_i(x,y,t)$, while the discrete
velocity vector field $\mathbf{v}_h$ is defined on the dual grid, namely  $\mathbf{v}_h(x,y,t)|_{\QQ_j^\st}=\mathbf{v}_j(x,y,t)$ . \par

The numerical solution of \eref{eq:CS_2}-\eref{eq:CS_2_2} is represented inside the space-time control volumes of the primal and the dual grid during the current time interval $T^{n+1}$ 
by piecewise space-time polynomials as follows: 
\begin{equation}
	p_i(x,y,t)=\sum\limits_{l=1}^{N_\phi^\st} \tilde{\phi}_l^{(i)}(x,y,t)\hat{p}_{l,i}^{n+1}=:\tilde{\bphi}^{(i)}(x,y,t)\etah_i^{n+1},
\label{eq:D_1}
\end{equation}
\begin{equation}
	\mathbf{v}_j(x,y,t)=\sum\limits_{l=1}^{N_\psi^\st} \tilde{\psi}_l^{(j)}(x,y,t) \hat{\mathbf{v}}_{l,j}^{n+1}=:\tilde{\bpsi}^{(j)}(x,y,t)\vvh_j^{n+1},
\label{eq:D_3}
\end{equation}
where the vectors of basis functions $\tilde{\bphi}(x,y,t)$ and $\tilde{\bpsi}(x,y,t)$ are generated via the mappings 
from $\tilde{\bphi}(\xi,\gamma, \tau)$ on $T_{std}\times [0,1]$ and $\bpsi(\xi,\gamma, \tau)$ on $R_{std} \times [0,1]$, respectively.

A weak formulation of the continuity equation \eref{eq:CS_2} is obtained by multiplying it by $\tilde{\bphi}$ and integrating over a control volume
$\TT_i^\st$, for every $k=1\ldots N_\phi^\st$. The resulting weak formulation reads
\begin{equation}
\int\limits_{\TT_i^\st}{\tilde{\phi}_k^{(i)} \, \nabla \cdot \mathbf{v} \, dx dy dt}=0.
\label{eq:CS_4}
\end{equation}
Similarly, multiplication of the momentum equation \eref{eq:CS_2_2} by $\tilde{\bpsi}$ and integrating over a control volume $\QQ_j^\st$ yields
\begin{equation}
\int\limits_{\QQ_j^\st}{\tilde{\psi}_k^{(j)}\left( \diff{\mathbf{v}}{t}+\nabla \cdot \mathbf{\TF} \right)  dx dy dt}+\int\limits_{\QQ_j^\st}{\tilde{\psi}_k^{(j)} \nabla \np \, dx dy dt}=\int\limits_{\QQ_j^\st}{\tilde{\psi}_k^{(j)} \, \mathbf{S} \, dx dy dt},
\label{eq:CS_5}
\end{equation}
for every $j=1 \ldots N_j$ and $k=1 \ldots N_\psi^\st$. Using integration by parts Eq. \eref{eq:CS_4} becomes 
\begin{equation}
\oint\limits_{\partial \TT_i^\st}{\tilde{\phi}_k^{(i)} \mathbf{v} \cdot \nv_{i} \, ds dt}-\int\limits_{\TT_i^\st}{\nabla \tilde{\phi}_k^{(i)} \cdot \mathbf{v} \, dx dy dt}  =0,
\label{eq:CS_6}
\end{equation}
where $\nv_{i}$ indicates the outward pointing unit normal vector.
Due to the discontinuity of $p_h$ and $\mathbf{v}_h$ at element boundaries, equations \eref{eq:CS_5} and \eref{eq:CS_6} have to be split. Note, however, that thanks to the use of a 
\textit{staggered grid} we do \textit{not} need a Riemann solver here, since all the quantities are readily defined where needed for the flux computation. In other words, the velocity
is continuous across the boundaries of the triangles on the main grid and the pressure is continuous across the boundaries of the dual quadrilateral grid. 

\begin{equation}
\sum\limits_{j \in S_i}\left( \int\limits_{\Gamma_j^\st}{\tilde{\phi}_k^{(i)} \mathbf{v}_j \cdot \nextud \, ds dt}-\int\limits_{\TT_{i,j}^\st}{\nabla \tilde{\phi}_k^{(i)} \cdot \mathbf{v}_j \, dx dy dt}  \right)=0,
\label{eq:CS_8}
\end{equation}
and
\begin{eqnarray}
\int\limits_{\QQ_j^\st}{\tilde{\psi}_k^{(j)}\left( \diff{\mathbf{v}_j}{t}+\nabla \cdot \mathbf{\TF}_j \right)  dx dy dt}
+\hspace{-3mm} \int\limits_{\TT_{\ell(j),j}^\st}{\tilde{\psi}_k^{(j)} \nabla \np_{\ell(j)} dx dy dt} \nonumber
+\hspace{-3mm} \int\limits_{\TT_{r(j),j}^\st}{\tilde{\psi}_k^{(j)} \nabla \np_{r(j)} \, dx dy dt} + \nonumber \\
\int\limits_{\Gamma_j^\st}{\tilde{\psi}_k^{(j)} \left(\np_{r(j)}-\np_{\ell(j)}\right) \nstd \, ds dt}=\int\limits_{\QQ_j^\st}{\tilde{\psi}_k^{(j)} \mathbf{S} dx dy dt},
\label{eq:CS_9}
\end{eqnarray}
where $\nextud=\nv_{i}|_{\Gamma_j^\st}$; $\TT_{i,j}^\st=\TT_{i,j} \times T^{n+1}$; and $\Gamma_j^\st=\Gamma_j \times T^{n+1}$.
Using definitions \eref{eq:D_1} and \eref{eq:D_3}, we conveniently rewrite the above equations as
\begin{eqnarray}
\sum\limits_{j \in S_i}\left(\int\limits_{\Gamma_j^\st}{\tilde{\phi}_k^{(i)}\tilde{\psi}_l^{(j)} \nextud ds dt}\, \mathbf{\hat{v}}_{l,j}^{n+1}-\int\limits_{\TT_{i,j}^\st}{\nabla \tilde{\phi}_k^{(i)}\tilde{\psi}_l^{(j)}dx dy dt} \, \mathbf{\hat{v}}_{l,j}^{n+1} \right)=0,
\label{eq:CS_10}
\end{eqnarray}
and
\begin{eqnarray}
\int\limits_{\QQ_j^\st}{\tilde{\psi}_k^{(j)} \diff{\mathbf{v}_j}{t} dx dy dt} +
\int\limits_{\QQ_j^\st}{\tilde{\psi}_k^{(j)} \nabla \cdot \mathbf{\TF} dx dy dt} \nonumber \\
+\int\limits_{\TT_{\ell(j),j}^\st}{\tilde{\psi}_k^{(j)} \nabla \tilde{\phi}_{l}^{(\ell(j))}  dx dy dt} \,  \, \hat \np_{l,\ell(j)}^{n+1}
+\int\limits_{\TT_{r(j),j}^\st}{   \tilde{\psi}_k^{(j)} \nabla \tilde{\phi}_{l}^{(r(j))}     dx dy dt} \,  \, \hat \np_{l,r(j)}^{n+1}    \nonumber \\
 +\int\limits_{\Gamma_j^\st}{\tilde{\psi}_k^{(j)} \tilde{\phi}_{l}^{(r(j))}    \nstd ds dt} \,  \hat \np_{l,r(j)}^{n+1}
 -\int\limits_{\Gamma_j^\st}{\tilde{\psi}_k^{(j)} \tilde{\phi}_{l}^{(\ell(j))} \nstd ds dt} \,  \hat \np_{l,\ell(j)}^{n+1}=\int\limits_{\QQ_j^\st}{\tilde{\psi}_k^{(j)} \mathbf{S} dx dy dt}, 
\label{eq:CS_11_1}
\end{eqnarray}
where we have used the standard summation convention for the repeated index $l$.
Integrating the first integral in \eref{eq:CS_11_1} by parts in time we obtain 
\begin{eqnarray}
\int\limits_{\QQ_j^\st}{\tilde{\psi}_k^{(j)} \diff{\mathbf{v}_j}{t} dx dy dt} &=& 
 \left[ \int\limits_{\QQ_j}{\tilde{\psi}_k^{(j)} \mathbf{v}_j dx dy }  \right]_{t=t^{n+1}}-\left[ \int\limits_{\QQ_j}{\tilde{\psi}_k^{(j)} \mathbf{v}_j dx dy }  \right]_{t=t^{n}} 
- \int\limits_{\QQ_j^\st}{\diff{\tilde{\psi}_k^{(j)}}{t} \mathbf{v}_j dx dy dt}  
\label{eq:CS_11_2}
\end{eqnarray}
In Eq. \eref{eq:CS_11_2} we can recognize the fluxes between the current space-time element $\QQ_j \times T^{n+1}$, the future space-time slab and the past space-time elements, as well as an internal space-time volume contribution that connects the layers inside the space-time element $\QQ_j^\st$ in an asymmetric way. Note how the asymmetry affects only the space-time volume contribution in \eref{eq:CS_11_2}. This is due to the nature of the the time derivative operator, which has a natural positive direction given by the causality principle in time. 
By substituting Eq. \eref{eq:CS_11_2} into \eref{eq:CS_11_1} we obtain the following weak formulation of the momentum equation in space-time:
\begin{eqnarray}
 \left( \left[ \int\limits_{\QQ_j}{\tilde{\psi}_k^{(j)} \tilde{\psi}_l^{(j)} dx dy }  \right]_{t=t^{n+1}} 
- \int\limits_{\QQ_j^\st}{\diff{\tilde{\psi}_k^{(j)}}{t} \tilde{\psi}_l^{(j)} dx dy dt} \right) \mathbf{\hat{v}}_{l,j}^{n+1}   & & \nonumber \\
+\int\limits_{\TT_{\ell(j),j}^\st}{\tilde{\psi}_k^{(j)} \nabla \tilde{\phi}_{l}^{(\ell(j))}  dx dy} \,  \, \hat \np_{l,\ell(j)}^{n+1}
+\int\limits_{\TT_{r(j),j}^\st}{   \tilde{\psi}_k^{(j)} \nabla \tilde{\phi}_{l}^{(r(j))}     dx dy} \,  \, \hat \np_{l,r(j)}^{n+1}   
 +\int\limits_{\Gamma_j^\st}{\tilde{\psi}_k^{(j)} \tilde{\phi}_{l}^{(r(j))}    \nstd ds} \,  \hat \np_{l,r(j)}^{n+1}
 -\int\limits_{\Gamma_j^\st}{\tilde{\psi}_k^{(j)} \tilde{\phi}_{l}^{(\ell(j))} \nstd ds} \,  \hat \np_{l,\ell(j)}^{n+1}  & & 
\nonumber \\ 
 =  \left[ \int\limits_{\QQ_j}{\tilde{\psi}_k^{(j)} \tilde{\psi}_l^{(j)} dx dy }  \right]_{t=t^{n}} \mathbf{\hat{v}}_{l,j}^{n}  
- \int\limits_{\QQ_j^\st}{\tilde{\psi}_k^{(j)} \nabla \cdot \mathbf{\TF} dx dy}
+ \int\limits_{\QQ_j^\st}{\tilde{\psi}_k^{(j)} \mathbf{S} \, dx dy dt}, & & \nonumber \\
\label{eq:CS_11}
\end{eqnarray}

For every $i$ and $j$, Eqs. \eref{eq:CS_10}-\eref{eq:CS_11} are written in a compact matrix form as
\begin{eqnarray}
    \sum\limits_{j \in S_i}\D_{i,j}\vvh_j^{n+1}=0 \label{eq:CS_12},
\end{eqnarray}
and
\begin{eqnarray}
    \left(\Mpsi_j^+ - \Mpsi_j^\circ \right) \vvh_j^{n+1} - \Mpsi_j^-\vvh_j^{n} + \Upsilon_j(\mathbf{v}, \nabla \mathbf{v}) + \RM_j \etah_{r(j)}^{n+1}- \LM_j \etah_{\ell(j)}^{n+1} =\Ss_j, \label{eq:CS_12_1}
\end{eqnarray}
respectively, where:
\begin{eqnarray}
	\Mpsi_j^+ &=& \int\limits_{\QQ_j}{\tilde{\psi}_k^{(j)}(x,y,t^{n+1})\tilde{\psi}_l^{(j)}(x,y,t^{n+1})  dx dy}, \label{eq:MD_2} \\
    \Mpsi_j^- &=& \int\limits_{\QQ_j}{\tilde{\psi}_k^{(j)}(x,y,t^{n+1})\tilde{\psi}_l^{(j)}(x,y,t^{n})  dx dy}, \label{eq:MD_2_1} \\
    \Mpsi_j^\circ &=& \int\limits_{\QQ_j^\st}{\diff{\tilde{\psi}_k^{(j)}}{t} \tilde{\psi}_l^{(j)} dx dy dt}, \label{eq:MD_2_2} \\
    \Upsilon_j &=& \int\limits_{\QQ_j^\st}{\tilde{\psi}_k^{(j)} \nabla \cdot \mathbf{\TF} dx dy dt}
\end{eqnarray}

\begin{equation}
	\D_{i,j}=\int\limits_{\Gamma_j^\st}{\tilde{\phi}_k^{(i)}\tilde{\psi}_l^{(j)}\nextud ds dt}-\int\limits_{\TT_{i,j}^\st}{\nabla \tilde{\phi}_k^{(i)}\tilde{\psi}_l^{(j)}dx dy dt},
\label{eq:MD_3}
\end{equation}

\begin{equation}
	\RM_{j}=\int\limits_{\Gamma_j^\st}{\tilde{\psi}_k^{(j)} \tilde{\phi}_{l}^{(r(j))}\nstd ds dt}+\int\limits_{\TT_{r(j),j}^\st}{\tilde{\psi}_k^{(j)} \nabla \tilde{\phi}_{l}^{(r(j))}  dx dy dt},
\label{eq:MD_4}
\end{equation}

\begin{equation}
	\LM_{j}=\int\limits_{\Gamma_j^\st}{\tilde{\psi}_k^{(j)} \tilde{\phi}_{l}^{(\ell(j))}\nstd ds dt}-\int\limits_{\TT_{\ell(j),j}^\st}{\tilde{\psi}_k^{(j)} \nabla \tilde{\phi}_{l}^{(\ell(j))}  dx dy dt},
\label{eq:MD_5}
\end{equation}

\begin{equation}
	\Ss_j=\int\limits_{\QQ_j^\st}{\tilde{\psi}_k^{(j)} \mathbf{S} dx dy dt}.
\label{eq:MD_5_2}
\end{equation}
Remark how $\Mpsi_j^\circ$ introduces, for polynomial degrees $p_\gamma>0$, an asymmetric contribution in time. 
The action of matrices $\LM$ and $\RM$ can be generalized by introducing a new matrix $\Q_{i,j}$, defined as
\begin{equation}
	\Q_{i,j}=\int\limits_{\TT_{i,j}^\st}{\tilde{\psi}_k^{(j)} \nabla \tilde{\phi}_{l}^{(i)}  dx dy dt}-\int\limits_{\Gamma_j^\st}{\tilde{\psi}_k^{(j)} \tilde{\phi}_{l}^{(i)}\sigma_{i,j} \nstd ds dt},
\label{eq:MD_6}
\end{equation}
where $\sigma_{i,j}$ is a sign function defined by
\begin{equation}
	\sigma_{i,j}=\frac{r(j)-2i+\ell(j)}{r(j)-\ell(j)}.
\label{eq:SD_1}
\end{equation}
In this way $\Q_{\ell(j),j}=-\LM_j$ and $\Q_{r(j),j}=\RM_j$, and then Eq. \eref{eq:CS_12_1} becomes in terms of $\Q$
\begin{equation}
	\left(\Mpsi_j^+ - \Mpsi_j^\circ \right) \vvh_j^{n+1} - \Mpsi_j^-\vvh_j^{n} + \Upsilon_j(\mathbf{v}, \nabla \mathbf{v})  + \Q_{r(j),j}  \etah_{r(j)}^{n+1} + \Q_{\ell(j),j} \etah_{\ell(j)}^{n+1} =\Ss_j,
\label{eq:CS_12_2}
\end{equation}
or, equivalently,
\begin{equation}
\left(\Mpsi_j^+ - \Mpsi_j^\circ \right) \vvh_j^{n+1} - \Mpsi_j^-\vvh_j^{n} + \Upsilon_j(\mathbf{v}, \nabla \mathbf{v})  + \Q_{i,j} \etah_{i}^{n+1} + \Q_{\p(i,j),j} \etah_{\p(i,j)}^{n+1} =\Ss_j.
\label{eq:CS_13}
\end{equation}

In order to further ease notation, we will use the abbreviation $\Mpsi_j = \Mpsi_j^+ - \Mpsi_j^\circ$ henceforth and will write  
Eqs. \eref{eq:CS_12}-\eref{eq:CS_12_1} as follows: 
\begin{eqnarray}
	\sum\limits_{j \in S_i}\D_{i,j}\vvh_j ^{{n+1}}=0,\label{eq:CS_15}	\\
	\Mpsi_j	\vvh_j^{n+1}-\Mpsi_j	\Fvh_j + \Q_{r(j),j} \etah_{r(j)}^{{n+1}}+ \Q_{\ell(j),j} \etah_{\ell(j)} ^{{ n+1}} =0,
\label{eq:CS_16}
\end{eqnarray}
where $\Fvh_j$ is an appropriate discretization of the nonlinear convective, viscous and source terms.
The details for the computation of $\Fvh_j$  will be presented later. 
Formal substitution of the discrete momentum equation \eref{eq:CS_16} into the discrete continuity equation \eref{eq:CS_15}, see also \cite{CasulliCheng1992,DumbserCasulli,2DSIUSW,2STINS}, yields 
\begin{eqnarray}
 \sum\limits_{j\in S_i}\D_{i,j}\Mpsi_j^{-1}\Q_{i,j} \etah_i^{n+1}
 +\sum\limits_{j\in S_i} \D_{i,j}\Mpsi_j ^{-1}\Q_{\p(i,j),j} \etah_{\p(i,j)}^{n+1}=\sum\limits_{j \in S_i} \D_{i,j} \Fvh_j,
\label{eq:CS_19}
\end{eqnarray}
We have now to choose a time discretization for the nonlinear convective-viscous term.
The simplest choice would be to take $\Fvh_j$ explicitly, so in this case $\sum\limits_{j \in S_i} \D_{i,j} \Fvh_j^{n}$ becomes a known term at time $t^n$ and hence Eq. \eref{eq:CS_19} would represent a four-point block system for the new pressure $\etah_i^{n+1}$, as proposed in \cite{2STINS}. Unfortunately, in problems when the convective-viscous effects cannot be neglected, this will produce only a low  order accurate method in time. The problem in this case is that the convective-viscous contribution in the time interval $T^{n+1}$ is based on the old information $T^n$ and does not see the effects of 
the new pressure in the time interval $T^{n+1}$. Furthermore, if we take $\Fvh_j$ implicitly, then system \eref{eq:CS_19} becomes nonlinear and it would be very cumbersome to solve it. In order to  overcome this problem we introduce a simple \textit{Picard iteration} to introduce the information of the new pressure into the viscous and convective terms, but without introducing a nonlinearity
in the final system to be solved. This approach is inspired by the local space-time Galerkin predictor method proposed for the high order time discretization of $P_NP_M$ schemes in \cite{Dumbser2008,ADERNSE}. 
Hence, for $k=1, N_{pic}$, we rewrite system \eref{eq:CS_19} as 
\begin{eqnarray}
 \sum\limits_{j\in S_i}\D_{i,j}\Mpsi_j^{-1}\Q_{i,j} \etah_i^{n+1,k+1}
 +\sum\limits_{j\in S_i} \D_{i,j}\Mpsi_j ^{-1}\Q_{\p(i,j),j} \etah_{\p(i,j)}^{n+1,k+1}=\sum\limits_{j \in S_i} \D_{i,j} \Fvh_j^{n+1,k \poh},
\label{eq:CS_19_2}
\end{eqnarray}
or, by introducing the boundary elements (see e.g. \cite{2STINS}),
\begin{eqnarray}
\left[ \sum\limits_{j\in S_i\cap \B(\Omega)}\D_{i,j}^{\partial}\Mpsi_j^{-1}\Q_{i,j}^{\partial} -\sum\limits_{j\in S_i-\B(\Omega)}\D_{i,j}\Mpsi_j^{-1}\Q_{i,j} \right] \etah_i^{n+1,k+1} 
- \sum\limits_{j\in S_i-\B(\Omega)} \D_{i,j}\Mpsi_j^{-1}\Q_{\p(i,j),j} \etah_{\p(i,j)}^{n+1,k+1}= \nonumber \\
-\sum\limits_{j \in S_i-\B(\Omega)} \D_{i,j} \Fvh_j^{n+1,k \poh}  +\sum\limits_{j \in S_i\cap\B(\Omega)} \D_{i,j}^{\partial} \Fvh_j^{n+1,k \poh},
\label{eq:82_2}
\end{eqnarray}

where $\D_{i,j}^{\partial}$ and $\Q_{i,j}^{\partial}$ are the natural extension of $\D$ and $\Q$ on triangular dual boundary elements, see e.g. \cite{2STINS}. 
Now the right hand side of Eq. \eref{eq:CS_19} can be computed by using the velocity field at the old Picard iteration $k$ and including the viscous effects using a fractional step type procedure.  
In this way, Eq. \eref{eq:CS_19} represents a block four-point system for the new pressure $\etah_i^{n+1,k+1}$. Once the new pressure field is known, the velocity vector field at the new 
Picard iteration $\vvh^{n+1,k+1}$ can be readily updated from the momentum equation \eref{eq:CS_16}.

\subsection{Nonlinear convection-diffusion}
To close the problem it remains to specify how to construct the nonlinear convection-diffusion operator $\Fvh_j^{n\poh}$. Following the ideas of \cite{2STINS}, a space-time DG scheme 
for the convection-diffusion terms on the dual mesh is given by 
\begin{eqnarray}
\int\limits_{\QQ_j^\st} \tilde{\psi}_k \frac{\partial}{\partial t} \mathbf{v}_h \, dx dy dt + \int\limits_{\partial \QQ_j^\st}{\tilde{\psi}_k \mathbf{G}_h \cdot \vec{n} \, \,  ds dt} - \int\limits_{\QQ_j^\st}{\nabla \tilde{\psi}_k \cdot \mathbf{F}(\mathbf{v}_h,\nabla \mathbf{v}_h) dx dy dt} 
 = \int\limits_{\QQ_j^\st}{\tilde{\psi}_k^{(j)} \mathbf{S} dx dy dt},
\label{eq:59}
\end{eqnarray}
\noindent and the numerical flux for both, the convective and the viscous contribution, is given such as in \cite{Rusanov:1961a,MunzDiffusionFlux,ADERNSE}, and reads
\begin{equation}
	\mathbf{G}_h \cdot \vec{n} = \frac{1}{2}\left(\mathbf{F}(\mathbf{v}_h^{\,+},\nabla \mathbf{v}_h^{\,+}) + \mathbf{F}(\mathbf{v}_h^{\,-},\nabla \mathbf{v}_h^{\,-}) \right)\cdot \vec{n} -\frac{1}{2}s_{\max} \left( \mathbf{v}_h^{\,+} - \mathbf{v}_h^{\,-} \right),
\label{eq:61}
\end{equation}

\noindent with

\begin{equation}
s_{\max} = 2 \, \max( |\mathbf{v}_h^{\,-} \cdot \vec n|, |\mathbf{v}_h^{\,+} \cdot \vec n| ) + \frac{2 \nu}{h^+ +h^-}\frac{2p+1}{\sqrt{\frac{\pi}{2}}},
\end{equation}
which contains the maximum eigenvalue of the Jacobian matrix of the purely convective transport operator $\mathbf{F}_c$ in normal direction,
see \cite{DumbserCasulli}, and the
stabilization term for the viscous flux, see \cite{ADERNSE,MunzDiffusionFlux}. Furthermore, the $\mathbf{v}_h^\pm$ and $\nabla \mathbf{v}_h^\pm$ denote the
velocity vectors and their gradients, extrapolated to the boundary of $\QQ_j$ from within the element $\QQ_j$ and from the neighbor element, respectively.
$h^+$ and $h^-$ are the maximum radii of the inscribed circle in $\QQ_j$ and the neighbor element, respectively.
We discretize the velocity $\mathbf{v}_h$ explicitly but its gradient has to be taken implicitly, in order to avoid additional restrictions on the maximum time step 
given by the viscous terms. In viscosity dominated problems, this allows us to use both, high viscosity and large time steps. After integration of the first term of \eqref{eq:59} 
by parts in time the resulting fully discrete formulation of \eref{eq:59} becomes  
\begin{equation}
 \hat{\mathbf{v}}_j^{n+1,k  \poh} = \Mpsi_j^{-1}\Mpsi_j^{-}\vvh_j^{n}- \Mpsi_j^{-1}\Upsilon_j(\mathbf{v}_h^{n+1,k},\nabla \mathbf{v}_h^{n+1,k  \poh})+\Mpsi_j^{-1}\Ss_j,
\label{eq:Ad2}
\end{equation}
where
\begin{eqnarray}
 \Upsilon_j(\mathbf{v}_h,\nabla \mathbf{v}_h) &=& \int\limits_{\QQ_j^\st}{\tilde{\psi}_k^{(j)} \nabla \cdot \mathbf{\TF}(\mathbf{v}_h,\nabla \mathbf{v}_h) dx dy} \nonumber \\
 	&=&  \int\limits_{\partial \QQ_j^\st}{\tilde{\psi}_k \mathbf{G}_h \cdot \vec{n} \, \,  ds} - \int\limits_{\QQ_j^\st}{\nabla \tilde{\psi}_k \cdot \mathbf{F}(\mathbf{v}_h,\nabla \mathbf{v}_h) dx dy}. 
\label{eq:Ad3}
\end{eqnarray}
Due to the explicit treatment of the nonlinear convective terms, the above method requires that the time step size is restricted by a CFL-type restriction for DG schemes, namely:
\begin{equation}
    \Delta t = \frac{\textnormal{CFL}}{2p+1}\cdot \frac{h_{min}}{2|\mathbf{v}_{max}|},
\label{eq:CFLC}
\end{equation}
where $h_{min}$ is the smallest incircle diameter; $\textnormal{CFL}<0.5$; and $\mathbf{v}_{max}$ is the maximum convective speed.
Furthermore, the time step of the global semi-implicit scheme is \textit{not} affected by the local time step used for the time integration of the
convective terms if a local time stepping / subcycling approach is employed, see \cite{CasulliZanolli,TavelliDumbserCasulli}. Implicit discretization
of the viscous contribution $\nabla \mathbf{v}$ in \eref{eq:59} involves two five-point block systems (one for each velocity component) that can be efficiently 
solved using a matrix-free GMRES algorithm \cite{GMRES}. The solution of this system is not necessary in problems where the viscous term is small enough to be integrated 
explicitly in time. In that case, i.e. for explicit discretizations of the viscous terms, one has to include the additional explicit time step restriction for 
parabolic PDE in eq. \eref{eq:CFLC}. 

Once $\mathbf{v}_j^{n+1,k \poh}$ has been computed, we set $\Fvh_j^{n+1,k  \poh} := \hat{\mathbf{v}}_j^{n+1,k  \poh}$.
As initial guess $\hat{\mathbf{v}}_j^{n+1,0}$ we can take the old velocity $\mathbf{v}_h^{n}$, or the extrapolation of $\mathbf{v}_h^{n}$ into the interval $T^{n+1}$. 

\subsection{Pressure correction formulation and final algorithm}
\label{sec_2.6}
The preliminary algorithm described above, as formulated by Eqs. \eref{eq:Ad2}, \eref{eq:82_2}, \eref{eq:CS_16} still contains an important drawback: 
indeed, Eq. \eref{eq:Ad2} does not depend on the pressure of the previous Picard iteration and hence the algorithm does not see the effect of the pressure in the time interval $T^{n+1}$. In order to overcome the problem we introduce  the contribution of the pressure from the previous Picard iteration directly into Eq. \eref{eq:Ad2}. Then, we update the velocity with the new pressure $\etah_i^{n+1,k+1}$. With this modification,  
Eqs. \eref{eq:Ad2},  \eref{eq:82_2}, \eref{eq:CS_16} and hence the final algorithm become: 

\begin{eqnarray}
\hat{\mathbf{v}}_j^{n+1,k  \poh} &=& \Mpsi_j^{-1}\Mpsi_j^{-}\vvh_j^{n}- \Mpsi_j^{-1}\Upsilon_j(\mathbf{v}_h^{n+1,k},\nabla \mathbf{v}_h^{n+1,k  \poh}) 
 - \Q_{r(j),j} \etah_{r(j)}^{{n+1,k}}- \Q_{\ell(j),j} \etah_{\ell(j)} ^{{ n+1,k}}+\Mpsi_j^{-1}\Ss_j, 
\label{eq:A1_mod}
\end{eqnarray}

\begin{eqnarray}
 \sum\limits_{j\in S_i}\D_{i,j}\Mpsi_j^{-1}\Q_{i,j} \left( \etah_i^{n+1,k+1} - \etah_i^{n+1,k} \right) 
 +\sum\limits_{j\in S_i} \D_{i,j}\Mpsi_j ^{-1}\Q_{\p(i,j),j} \left( \etah_{\p(i,j)}^{n+1,k+1} - \etah_{\p(i,j)}^{n+1,k} \right) 	
 =\sum\limits_{j \in S_i} \D_{i,j}  \Fvh_j^{n+1,k \poh}, 
\label{eq:A2_mod}
\end{eqnarray}

\begin{eqnarray}
		\vvh_j^{n+1,k+1}&=&\Fvh_j^{n+1,k \poh} 
		- \Mpsi_j^{-1}\left[\Q_{r(j),j} \left(\etah_{r(j)}^{{n+1,k+1}}- \etah_{r(j)}^{{n+1,k}} \right)- \Q_{\ell(j),j} \left(\etah_{\ell(j)} ^{{ n+1, k+1}}-\etah_{\ell(j)} ^{{ n+1, k}}\right)\right].
\label{eq:A3_mod}
\end{eqnarray}
Note that Eqs. \eqref{eq:A2_mod} and \eqref{eq:A3_mod} are written in terms of the \textit{pressure correction} 
$\Delta \etah_i^{{ n+1, k+1}} = \left(\etah_i^{{ n+1, k+1}}-\etah_i^{{ n+1, k}}\right)$.

As initial guess for the pressure one can take $\mathbf{p}_h^{{n+1,0}}=0$, but one could also choose the extrapolation of $\mathbf{p}_h^{n}$ into $T^{n+1}$.
One time step of the final algorithm can be summarized as follows: 
\begin{enumerate}
	\item Initialize $\mathbf{v}_h^{n+1,0}$ and $\mathbf{p}_h^{{n+1,0}}$ using the known information from $T^n$; 
	\item Picard iteration over $k=0\ldots N_{pic}$:
		\begin{enumerate}
			\item compute $\mathbf{v}_h^{n+1,k  \poh}$ using \eref{eq:A1_mod}, i.e. convective terms are discretized explicitly and viscous terms implicitly; then 
			set $\Fvh_j^{n+1,k \poh}:=\hat{\mathbf{v}}_j^{n+1,k \poh}$, 
			\item compute $\etah^{{n+1,k+1}}$ by solving the discrete pressure Poisson equation \eref{eq:A2_mod},  
			\item update $\vvh_j^{n+1,k+1}$ explicitly from \eref{eq:A3_mod};  
		\end{enumerate}	
	\item set $\vvh_j^{n+1}=\vvh_j^{n+1,k+1}$	and  $\etah^{{n+1}}=\etah^{{n+1,k+1}}$. 
\end{enumerate}

For the spatial computational domain we can apply the remark given in \cite{2STINS} and so either use a subparametric or a complete isoparametric representation. The second approach requires to store 
more information about each element, but it also allows to generalize the shape of the elements. This property is crucial when we try to discretize complex curved domains with a very coarse grid. 
In any case, this generalization does not affect the computational time during run-time, since it interests only the construction of the geometry-dependent matrices in the preprocessing stage of the
algorithm. 

\subsection{Splitting of the space-time matrices into a spatial and temporal part}
Even if the shape of the main matrices is similar compared to the ones introduced in \cite{2STINS}, the number of degree of freedom and the integral values are, in general, different. Due to the tensor  product construction of the space-time basis functions, we can split the main integrals \eref{eq:MD_2}-\eref{eq:MD_3} and \eref{eq:MD_6} into a spatial and a temporal part. Briefly, the space-time  matrices are generated from the spatial matrices of \cite{2STINS}, componentwise, as: 

\begin{eqnarray}
	\Mpsi_j^+(k,l) &=& \gamma_{\ell_2(k)}(t^{n+1})\gamma_{\ell_2(k)}(t^{n+1})\Mpsi_j^s(\ell_1(k),\ell_1(l)), \\
    \Mpsi_j^-(k,l) &=& \gamma_{\ell_2(k)}(t^{n+1})\gamma_{\ell_2(k)}(t^{n})\Mpsi_j^s(\ell_1(k),\ell_1(l)), \\
    \Mpsi_j^\circ(k,l)  &=& \Mpsi_j^s(\ell_1(k),\ell_1(l))\D^t(\ell_2(k), \ell_2(l)), \\
    \D_{i,j}(k,l)&=& \Delta t^{n+1} \D_{i,j}^s(\ell_1(k), \ell_1(l))\Mpsi^t(\ell_2(k), \ell_2(l)), \\
    \Q_{i,j}(k,l) &=& \Delta t^{n+1} \Q_{i,j}^s(\ell_1(k), \ell_1(l))\Mpsi^t(\ell_2(k), \ell_2(l)),
\end{eqnarray}
where the apex $s$ means that the matrix is the one constructed in \cite{2STINS}; $\D^t$ and $\Mpsi^t$ are two time matrices defined as 
\begin{eqnarray}
	 	\D^t\left(\tilde{k}, \tilde{l} \right) &=& \int\limits_0^1 \frac{d\gamma_{\tilde{k}}(\xi)}{d\xi}\gamma_{\tilde{l}}(\xi) d\xi , \label{eq:Dt}\\
		\Mpsi^t\left(\tilde{k}, \tilde{l} \right)&=&\int\limits_0^1 \gamma_{\tilde{k}}(\xi)\gamma_{\tilde{l}}(\xi) d\xi,
\end{eqnarray}
Remark how the action of the matrix $\D^t$ defined in \eref{eq:Dt} is symmetric only if $p_\gamma=0$.

\newpage 

\section{Numerical test problems}
In this section we study the accuracy of our new numerical method by solving some classical numerical benchmark problems, such as the lid-driven cavity flow, the unsteady oscillatory flow in a pipe or the unsteady flow past a circular cylinder. In particular, we perform quantitative comparisons between the numerical solution and available exact analytical solutions wherever possible. 
\label{sec.tests}

\subsection{Convergence test using a manufactured solution}

In order to study the accuracy of the proposed space-time DG method, we need an exact unsteady solution of \eref{eq:CS_2_2_0}-\eref{eq:CS_2_2}. For that purpose, we propose a so-called 
\textit{manufactured solution} in this section, which also makes use of a linear source term of the type $\mathbf{S}(x,y,t)$. The exact analytical solution for the velocity and the pressure 
is constructed so that 
\begin{equation}
    \mathbf{v}_{an} = \mathbf{v}_0 \sin\left[k(x-y)-\omega t \right], \qquad 
    p_{an} = p_0 \sin\left[k(x-y)-\omega t \right],
    \label{eq:NT_1_1}
\end{equation}
with the amplitudes $\mathbf{v}_0=(u_0,v_0)$ and $p_0$. Using the manufactured solution $(\mathbf{v}_{an},p_{an})$ we can compute all terms in \eref{eq:CS_2_2_0} exactly and hence 
obtain a source term $\mathbf{S}(x,y,t)$ that balances the momentum equation. 
Remark that the velocity field must be divergence-free ($\nabla \cdot \mathbf{v}=0$), hence $u_0=v_0$. In the present test case, we take $u_0=v_0=1$; $p_0=1$; $\omega=2 \pi$; $k=10/2\pi$; $t_{end}=0.5$; $\Delta t$ according to condition \eref{eq:CFLC}; and $\nu=0.01$. The temporal accuracy is chosen equal to the spatial one, the total number of Picard iterations is taken as $N_{pic}=p+1$ and $p^{n+1,0}  \equiv 0$ for the present test.  
The computational domain is $\Omega=[-0.5 , 0.5]^2$; the exact velocity field and pressure are taken as initial conditions and the exact pressure is also specified on $\partial \Omega$ 
as boundary condition. The $L_2$ error between the analytical and the numerical solution is computed as 
$$
\epsilon(p)=\sqrt{\int_\Omega {(p_{h}-p_{an})^2 dx dy}} \quad , \quad \epsilon(\mathbf{v})=\sqrt{\int_\Omega {(\mathbf{v}_{h}-\mathbf{v}_{an})^2 dx dy}}
$$
where the subscript $h$ refers to the numerical solution obtained at the final time $t=t_{end}$. The resulting rate of convergence is shown in Table \eref{tab:2}. 
We observe that the optimal order of convergence is obtained up to $p=4$ for the present unsteady test. 
 
\begin{table}[!htb]
\caption{Numerical convergence results for the manufactured solution test problem with polynomial degrees $p=1$ to $p=4$ in space and time.}
\begin{center}
\begin{tabular}{ccccccccc}
\hline
$N_i$ & \multicolumn{1}{c}{$\epsilon(p)$} & \multicolumn{1}{c}{$\epsilon(\mathbf{v})$} & \multicolumn{1}{c}{$\sigma (p)$} & \multicolumn{1}{c}{$\sigma (\mathbf{v})$} & \multicolumn{1}{c}{$\epsilon(p)$} & \multicolumn{1}{c}{$\epsilon(\mathbf{v})$} & \multicolumn{1}{c}{$\sigma (p)$} & \multicolumn{1}{c}{$\sigma (\mathbf{v})$} \\
\hline  
 & \multicolumn{4}{c}{$p=p_\gamma=1$} & \multicolumn{4}{c}{$p=p_\gamma=2$} \\  
  \hline  
    40 &  1.217E-01 & 9.572E-02 & - & -        &  8.740E-03   &  1.052E-02         &  -     &   -    \\
    160 &  2.678E-02 & 2.362E-02  & 2.2 & 2.0   &  8.833E-04   & 1.065E-03     & 3.3  &  3.3 \\
    640 &  6.050E-03 & 5.527E-03 & 2.1 & 2.1   &  1.050E-04   & 9.103E-05      & 3.1  &  3.5 \\
    2560 &  1.758E-03 & 1.497E-03 & 1.8 & 1.9   &  1.347E-05   & 7.820E-06      & 3.0  &  3.5 \\
  \hline
   & \multicolumn{4}{c}{$p=p_\gamma=3$} & \multicolumn{4}{c}{$p=p_\gamma=4$} \\
  \hline  
    40 &  7.703E-04 & 1.425E-03  & - & -        &  5.315E-05   &  7.135E-05         &  -     &   -    \\
    160 & 3.864E-05& 4.999E-05  & 4.3 & 4.8   &  1.143E-06   & 1.418E-06     & 5.5  &  5.7 \\
    640 & 2.425E-06 & 1.974E-06  & 4.0 & 4.7   &  3.102E-08   & 2.945E-08      & 5.2  &  5.6 \\
  \hline
\end{tabular}
\end{center}
\label{tab:2}
\end{table}

\subsection{The Womersley problem}
Here we consider an unsteady, viscosity-dominated test problem for which the incompressible Navier-Stokes equations have a nontrivial exact solution, namely the fluid flow 
inside a rigid planar pipe that is driven by a sinusoidal pressure gradient of the type 
\begin{equation}
    \frac{p_{out}(t)-p_{in}(t)}{L}=\Re \left(\frac{\tilde{P}}{\rho} e^{i\omega t}\right). 
    \label{eq:W_1}
\end{equation}
In this test $L$ denotes the tube length; $\tilde{P}$ is the amplitude of the pressure oscillation; $\rho$ is the density of the fluid; $\omega$ is the frequency of the oscillation; $p_{in}$ and $p_{out}$ indicate  the inlet and the outlet pressure, respectively; $\Re$ is the real part operator.
By imposing Eq. \eref{eq:W_1} at the tube ends, the exact analytical solution for the three dimensional, axially symmetric case was found by Womersley in \cite{Womersley1995}. It can be derived also for  the two dimensional planar case. The resulting axial velocity is uniform in the $x-$direction and is given by 
\begin{equation}
    u(x,y,t)=\Re \left[ i \frac{\tilde{P}}{\rho}{\omega} \left( 1- \frac{\cos[\lambda (y_c-1)]}{\cos(\lambda)} \right) \right],
    \label{eq:W_2}
\end{equation}
where $\lambda= \sqrt{-i \alpha^2}$; $\alpha=R\sqrt{\frac{\omega}{\nu}}$; $y_c=\frac{y-y_{b}}{R}$; and $y_b$ is the $y$ value of the bottom.

For the present test $\Omega=[-0.5, 0.5]\times [-0.2, 0.2]$; and $\frac{\tilde{P}}{\rho}=1$. We take a set of successively refined grids in order to show the convergence behaviour to the exact solution  with respect to the order $p$ in space and $p_\gamma$ in time. According to \cite{Womersley1995} the nonlinear convection effect is neglected for the present test. Thus, the stability of our scheme 
is not restricted by the CFL condition on the fluid velocity. Since we use very large time steps and a high viscosity coefficient in this test, the implicit treatment 
of the viscous terms is necessary to allow large time steps. In particular we choose $\nu=5 \cdot 10^{-2}$ and $t_{end}=1.5$. On the coarsest grid we use $\Delta t = t_{end}/6$, then the time step 
is reduced proportional to the spatial grid size. No-slip boundary conditions are imposed on the top and the bottom boundary, while the pressure \eref{eq:W_1} is imposed at the inlet and the outlet boundary
on the left and on the right, respectively. The number of Picard iterations is given by $N_p=p+1$ for all simulations. 

\begin{figure}[ht]
    \begin{center}
   \includegraphics[width=0.49\textwidth]{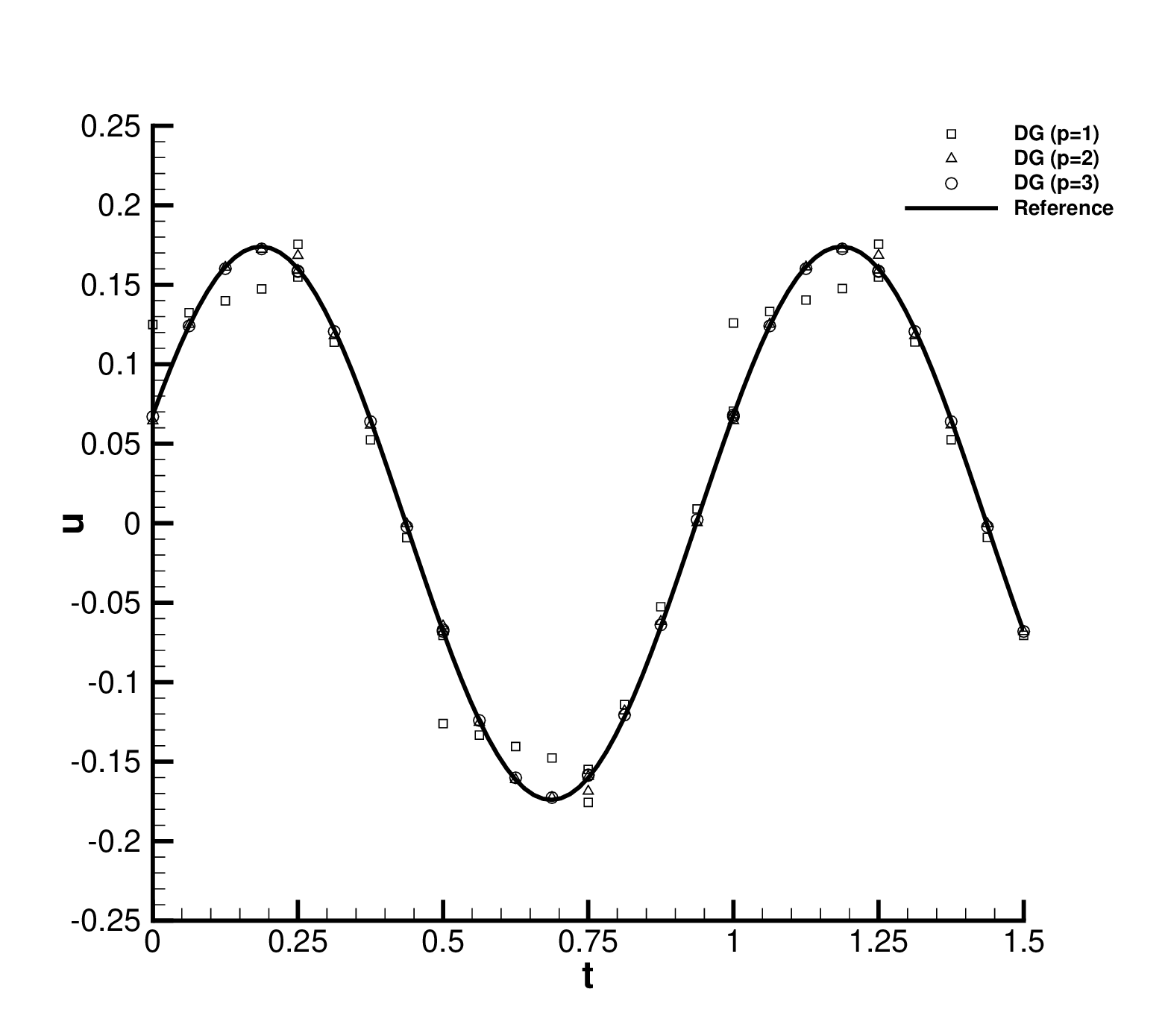}
   \includegraphics[width=0.49\textwidth]{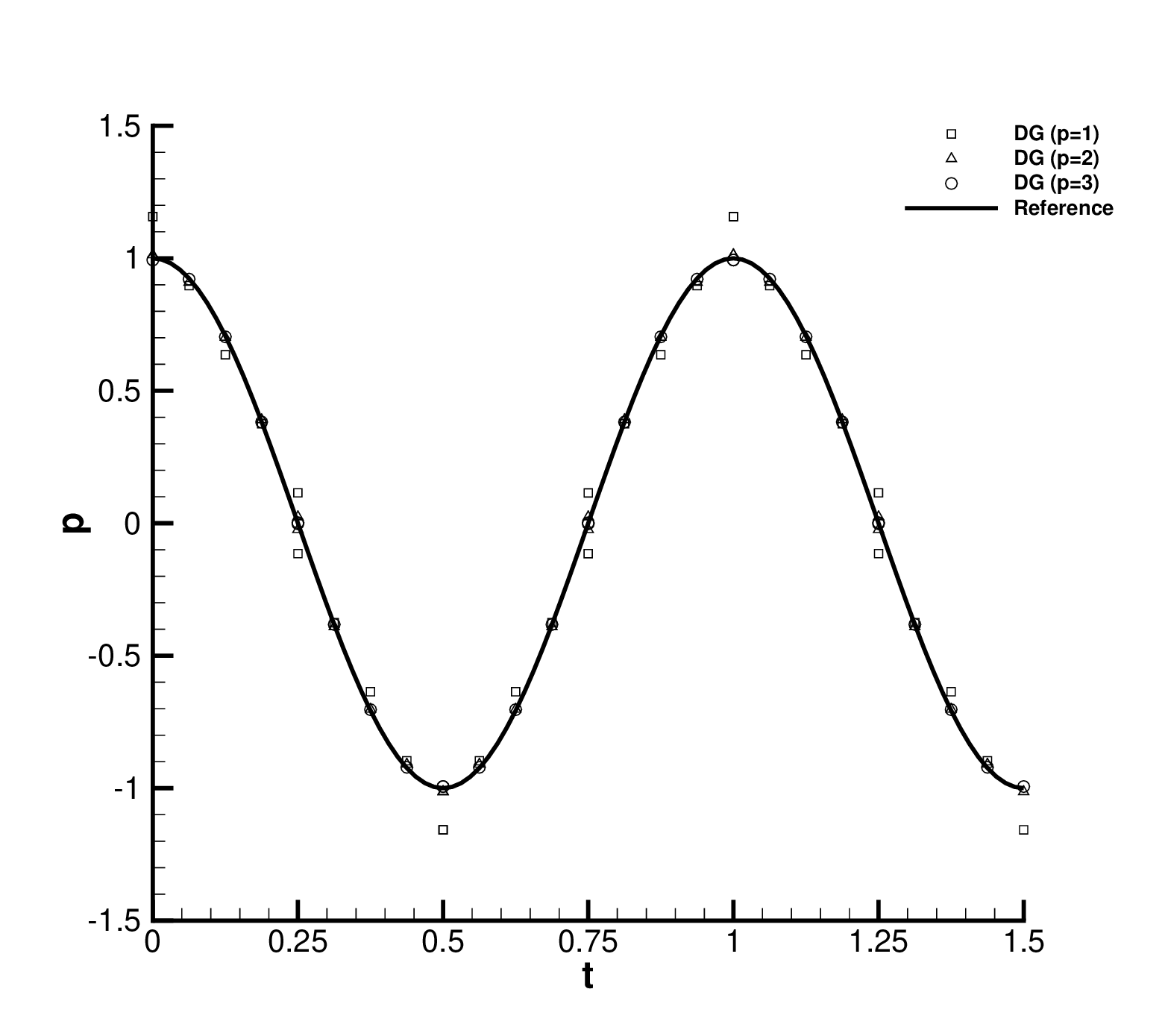}
    \caption{Time series for the axial velocity $u$ and the pressure $p$ computed at $(x,y)=(-0.5,0)$ for the coarsest grid $N_i=46$ and $N_t=6$}
    \label{fig.w1}
		\end{center}
\end{figure}

\begin{figure}[ht]
    \begin{center}
   \includegraphics[width=0.6\textwidth]{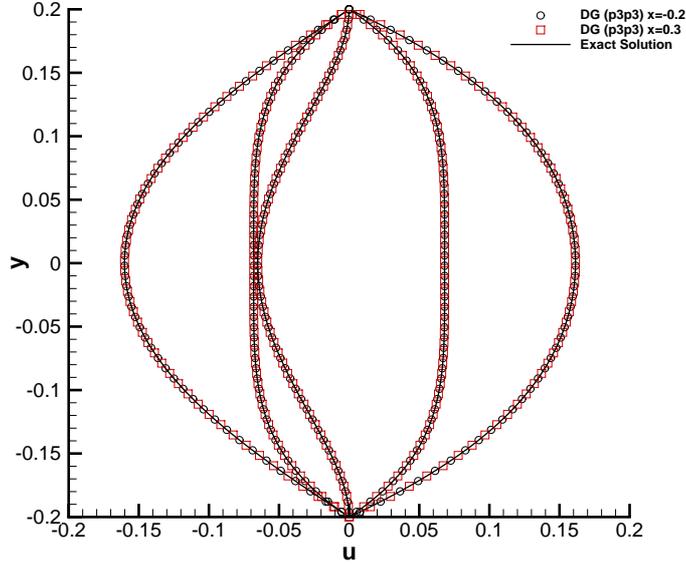}
    \caption{Radial velocity profiles for $x=-0.2$ and $x=0.3$ at times, 
from left to right, $t=[0.75,0.5,0.875,1.0,0.125]$. Comparison between exact and numerical solution.}
    \label{fig.w_n1}
                \end{center}
\end{figure}

\begin{table}[!htb]
\caption{Numerical convergence results for the planar Womersley problem.}
\begin{center}
\begin{tabular}{cccccccc}
	\hline
	$p$ & $p_{\gamma}$ & $N_i$ & $N_t$ & $\epsilon(p)$& $\epsilon(\mathbf{v})$ & $\sigma(p)$& $\sigma(\mathbf{v})$ \\
	\hline \hline
	1 & 1 & 46 & 6 & 5.7880182E-02&    1.8848423E-03  		& - & - 		  \\
	1 & 1 & 184 & 12 & 1.7635947E-02  & 5.5901107E-04 	& 1.7 & 1.8  \\
	1 & 1 & 736 & 24 & 4.6206559E-03   & 1.4587701E-04 & 1.9 & 1.9  \\
	1 & 1 & 2944 & 48 & 1.1683966E-03 & 3.7404869E-05 & 2.0 & 2.0 \\
    \hline
	2 & 2 & 46 & 6 & 7.0716231E-03&    2.6412698E-04 		& - & - 		  \\
	2 & 2 & 184 & 12 & 4.8160864E-04  & 3.8846170E-05	& 3.9 & 2.8  \\
	2 & 2 & 736 & 24 & 3.0677533E-05   & 7.2036760E-06 & 4.0 & 2.4  \\
	2 & 2 & 2944 & 48 & 1.9295385E-06 & 1.6070616E-06 & 4.0 & 2.2 \\
    \hline
	3 & 3 & 46 & 6 & 9.8372146E-04&   1.2793693E-05		& - & - 		  \\
	3 & 3 & 184 & 12 & 7.7144497E-05   & 7.8462176E-07 	& 3.7 & 4.0  \\
	3 & 3 & 736 & 24 & 5.0814347E-06   & 4.8795894E-08  & 3.9 & 4.0  \\
	3 & 3 & 2944 & 48 & 3.2173776E-07 & 3.0326872E-09 & 4.0 &4.0 \\
    \hline
	4 & 4 & 46 & 6 & 7.3692980E-05&   5.1193160E-07		& - & - 		  \\
	4 & 4 & 184 & 12 & 1.2539784E-06   & 2.1649081E-08 	& 5.9 & 4.6  \\
	4 & 4 & 736 & 24 & 2.1930727E-08   & 1.1576584E-09  & 5.8 & 4.2  \\
	4 & 4 & 2944 & 48 & 1.0258845E-09 & 7.0131498E-11 & 4.4 &4.0 \\
	\hline
\end{tabular}
\end{center}
\label{tab:3}
\end{table}

The resulting convergence results, using the $L_2-$norm as in the previous example, are shown in Table \ref{tab:3}. Observe how a non-optimal order of convergence $p$ is achieved for the velocity
for odd order schemes, while the optimal convergence rate $p+1$ is achieved for the pressure for all polynomial degrees. Note that when using the semi-implicit staggered DG method introduced 
in \cite{2DSIUSW} only a second order of convergence could be achieved for this unsteady test problem, while full high order convergence in space and time is obtained with the new scheme presented 
in this paper. 
In Figure \ref{fig.w1} we show the time series of the axial velocity and the pressure in a given point for the coarsest grid configuration $(N_i, N_t)=(46,6)$. While piecewise linear space-time 
polynomials are not able to reproduce the sinusoidal signal well with only six time steps, the piecewise quadratic and higher order approximations in space and time yield an almost perfect match 
with the exact solution even on this extremely coarse space-time grid. 

In Figure $\ref{fig.w_n1}$ we compare the resulting numerical velocity profiles $u(y)$ against the exact solution at several times for the case  
$(p,p_\gamma)=(3,3)$ and $N_i=736$. Two different locations, $x=-0.2$ and $x=0.3$, are plotted in order to show that the profile is constant in the  
$x$-direction. One observes that there is no visible difference between numerical and exact solution in Fig. \ref{fig.w_n1}.

\subsection{Taylor-Green vortex}
Another widely used testcase for the verification of numerical methods for the incompressible Navier-Stokes equations is the Taylor-Green vortex problem. 
The analytical unsteady solution is given by 
\begin{eqnarray}
    u(x,y,t)&=&\sin(x)\cos(y)e^{-2\nu t}, \label{eq:TG_0} \\
    v(x,y,t)&=&-\cos(x)\sin(y)e^{-2\nu t}, \label{eq:TG_1} \\
    p(x,y,t)&=&\frac{1}{4}(\cos(2x)+\cos(2y))e^{-4\nu t}.
\label{eq:TG_2}
\end{eqnarray}
The computational domain is $\Omega=[0,2\pi]^2$ and is extended using periodic boundary conditions on all the boundaries.

\begin{figure}[ht]
    \begin{center}
   \includegraphics[width=0.7\textwidth]{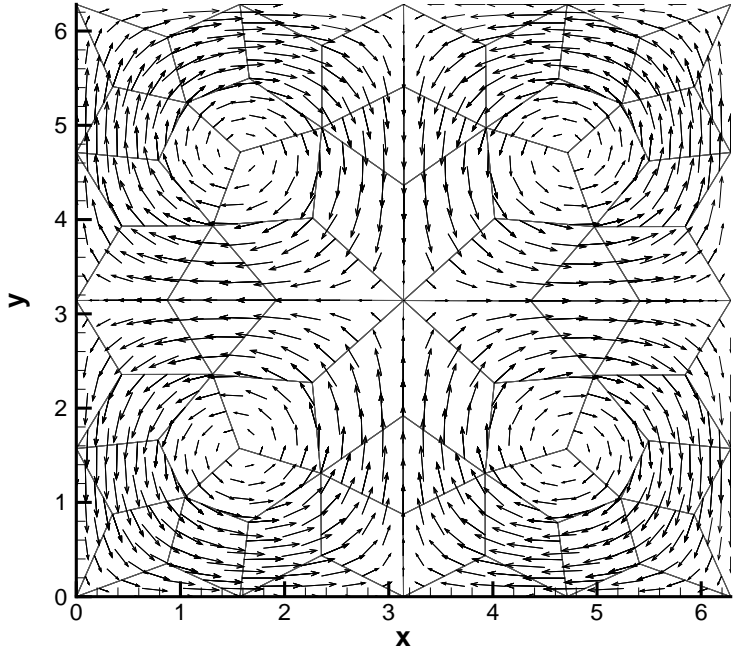}
    \caption{Velocity field of the Taylor-Green vortex on the coarse grid $N_i=40$ with $p=4$. The edge-based dual grid is shown.}
    \label{fig.tg1}
		\end{center}
\end{figure}
As implied by Eqs. \eref{eq:TG_0}-\eref{eq:TG_2}, the resulting velocity field initially appears as depicted in Figure \ref{fig.tg1} and then starts to lose energy according to the friction effects. 
For the present test we consider several grid refinements; $t_{end}=0.1$; $\nu=0.1$; and $\Delta t$ is chosen according to the CFL time restriction for the nonlinear convective terms. 
The numerical convergence results are shown in Table \ref{tab:4}. We find that the optimal convergence rates are achieved for this important nontrivial test problem with periodic boundary conditions. 

\begin{table}[!htb]
\begin{center}
\begin{tabular}{|c|c|c|c|c|c|c|c|c|}
\hline
$N_i$ & \multicolumn{2}{|c|}{$p=p_\gamma=1$} & \multicolumn{2}{c|}{$p=p_\gamma=2$}& \multicolumn{2}{|c|}{$p=p_\gamma=3$} & \multicolumn{2}{c|}{$p=p_\gamma=4$} \\
  & \multicolumn{1}{c}{$\epsilon(\mathbf{v})$} & $\sigma$ & \multicolumn{1}{c}{$\epsilon(\mathbf{v})$} & $\sigma$ & \multicolumn{1}{c}{$\epsilon(\mathbf{v})$} & $\sigma$ & \multicolumn{1}{c}{$\epsilon(\mathbf{v})$} & $\sigma$ \\
  \hline \hline
    40 	& 3.088E-01  	& --   	&  5.588E-02	& -- 	& 5.895E-03  & --  &  1.669E-03	& --      \\
    160 & 8.868E-02  	& 1.8   &  5.765E-03  & 3.3 & 4.730E-04  & 3.6 &  3.109E-05	& 5.7  \\
    640 & 2.267E-02  	& 2.0  	&  7.052E-04  & 3.0 & 2.387E-05  & 4.3 &  6.233E-07	& 5.6   \\
    2560& 5.476E-03 	& 2.0   &  8.452E-05	&	3.1	& 1.312E-06  & 4.2  &  1.297E-08           	& 5.6  \\
  \hline
\end{tabular}
\end{center}
\caption{Numerical convergence results for the velocity vector field of the Taylor-Green vortex.}
\label{tab:4}
\end{table}

\subsection{Double shear layer}
\begin{figure}[ht]
    \begin{center}
   \includegraphics[width=0.48\textwidth]{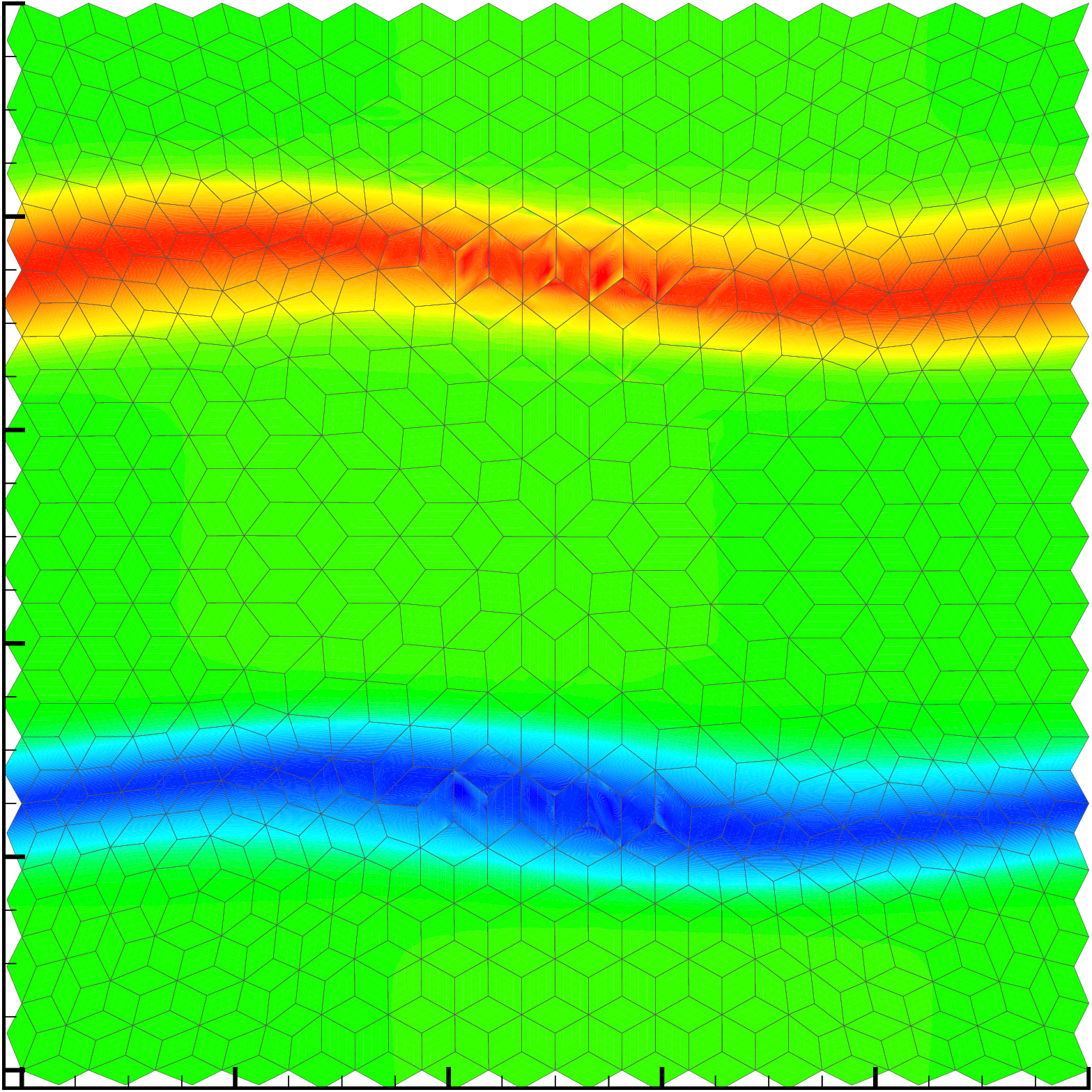}
   \includegraphics[width=0.48\textwidth]{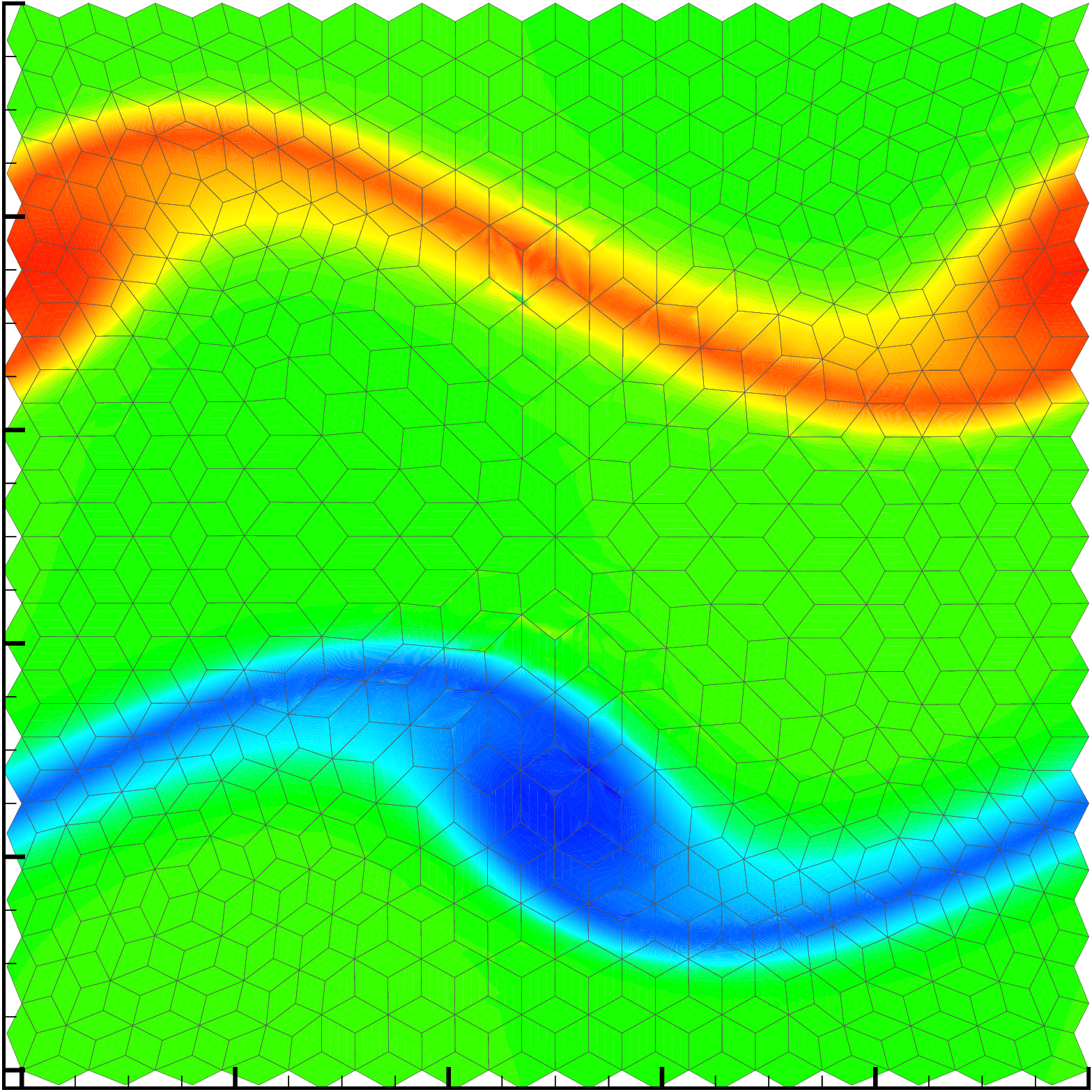}
   \includegraphics[width=0.48\textwidth]{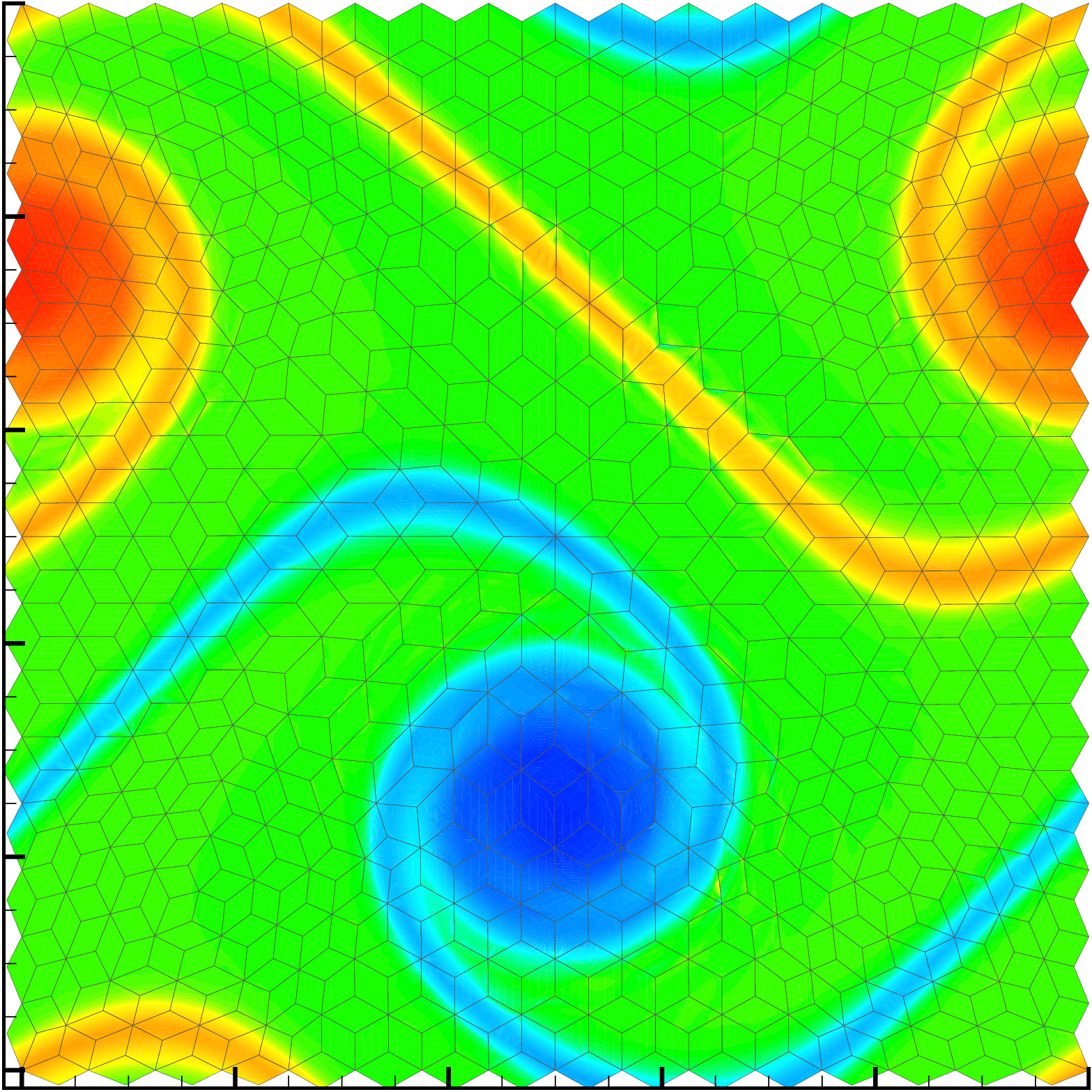}
   \includegraphics[width=0.48\textwidth]{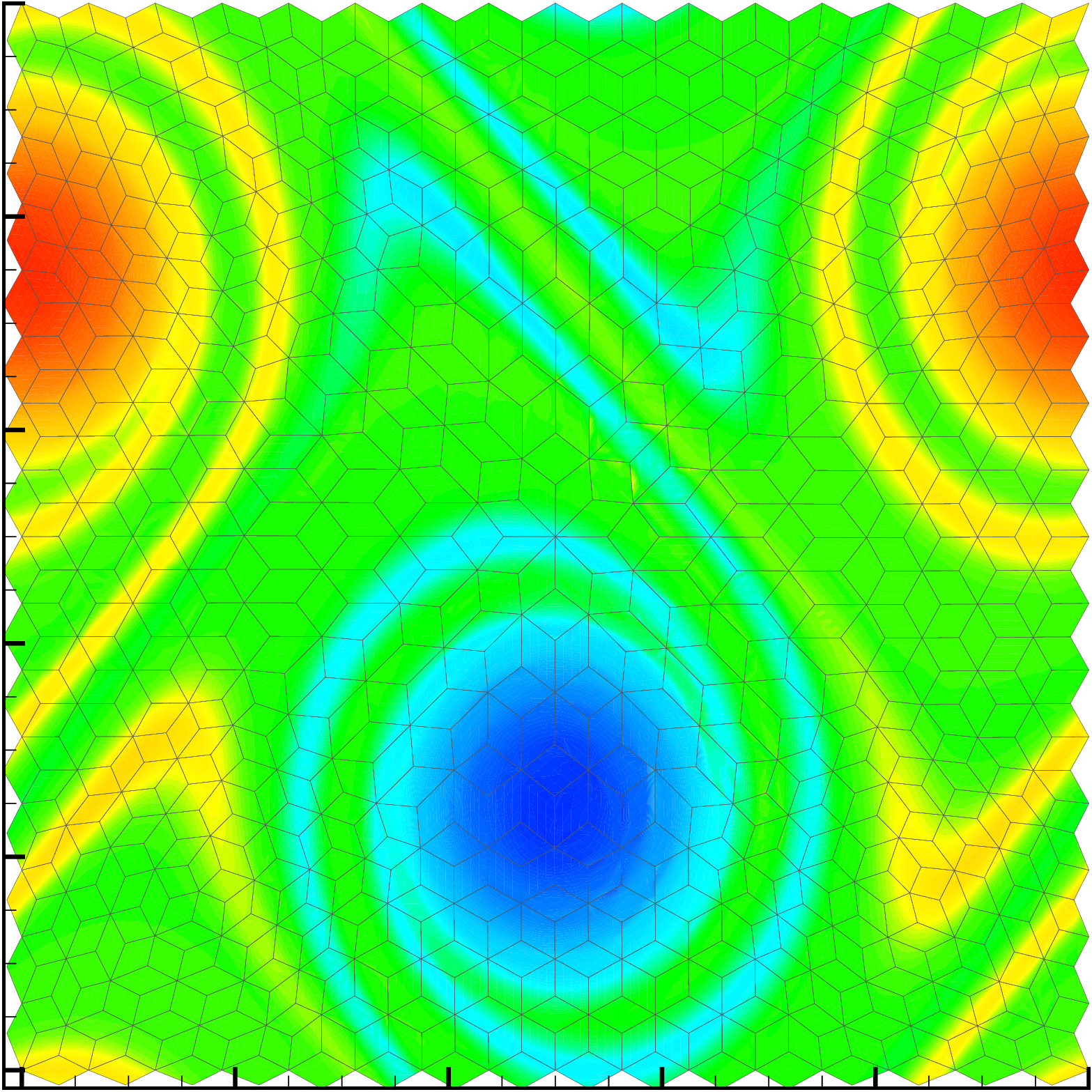}
    \caption{Vorticity pattern for the double shear layer test at times, from top left to bottom right, $t=0.4$;$t=0.8$;$t=1.2$;$t=1.8$}
    \label{fig.dsl1}
	\end{center}
\end{figure}

The numerical scheme is applied here to a test case studied in \cite{Bell1989}, which contains a high initial velocity gradient. 
We take $\Omega=[-1,1]^2$ and, as initial condition, we consider a perturbed double shear layer profile: 
\begin{eqnarray}
    u_0&=&\left\{
    \begin{array}{l}
    	\tanh\left[ \tilde{\rho} (y_n-0.25) \right], \qquad \textnormal{ if } y_n \leq 0.5, \\
    	\tanh\left[ \tilde{\rho} (0.75-y_n) \right], \qquad \textnormal{ if } y_n > 0.5,
    \end{array}
    \right.
    \label{eq:DSL_0} \\
    v_0&=& \delta \sin(2\pi x_n), \\
    p_0&=&1,
\label{eq:DSL_2}
\end{eqnarray}
where $y_n=\frac{y+1}{2}$ and $x_n=\frac{x+1}{2}$ are the normalized vertical and horizontal coordinates, respectively; $\tilde{\rho}$ is a parameter that determines the slope of the shear layer; and $\delta$ is the amplitude of the initial perturbation. For the present test we set $\delta=0.05$; $\tilde{\rho}=30$; $\nu=2\cdot10^{-4}$; $p=4$ and $p_\gamma=3$. The time step is chosen according to the CFL condition for the nonlinear convective terms and four Picard iterations have been used in this simulation. The domain $\Omega$ is covered with a total number of only $N_i=640$ triangles and periodic  boundary conditions are imposed everywhere. The resulting vorticity pattern is reported at several times in Figure \ref{fig.dsl1}. 
The two thin shear layers evolve into several vortices, as observed in \cite{Bell1989}, and overall the small flow structures seem to be relatively well resolved also at the final time $t=1.8$, 
even if a very coarse grid has been used in space and time.

\subsection{Lid-driven cavity flow}
We consider here another classical benchmark problem for the incompressible Navier-Stokes equations, namely the lid-driven cavity problem \cite{Ghia1982}. 
This test case is solved numerically with the new staggered space-time DG scheme on very coarse grids using polynomial degrees of $p=3$ 
and $p_{\gamma}=3$ in space and time, respectively. 
Let $\Omega=[-0.5,0.5]\times [-0.5, 0.5]$, set velocity boundary conditions $u=1$ and $v=0$ on the top boundary (i.e. at $y=0.5$) and 
impose so-slip wall boundary conditions on the other edges. As initial condition we take $u(x,y,0)=v(x,y,0)=0$. We use a grid with only 
$N_i=116$ triangles for $Re=100,400,1000$ and $N_i=512$ triangles for $Re=3200$.

\begin{figure}[ht]
    \begin{center}
    \includegraphics[width=0.48\textwidth]{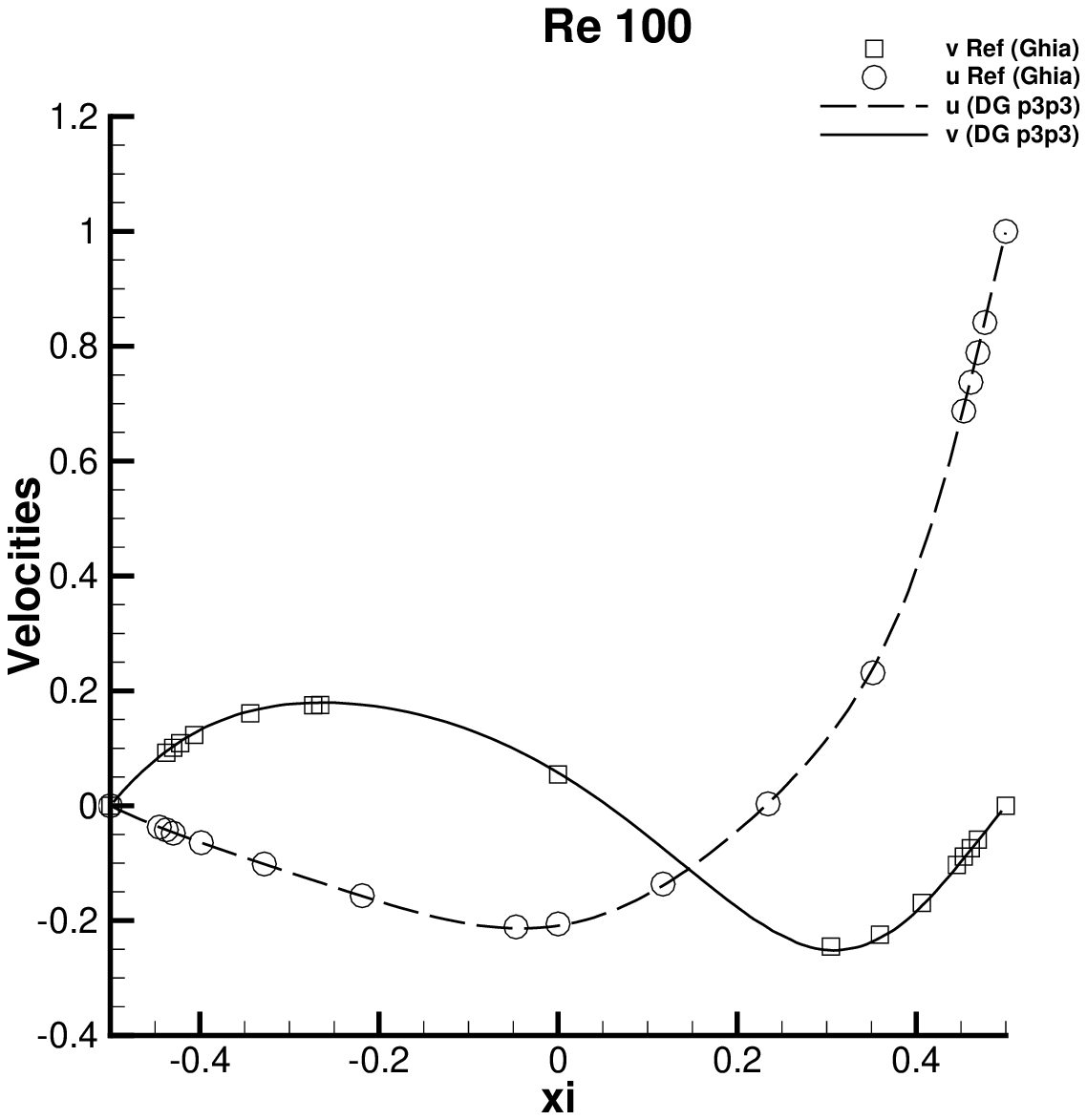}
    \includegraphics[width=0.48\textwidth]{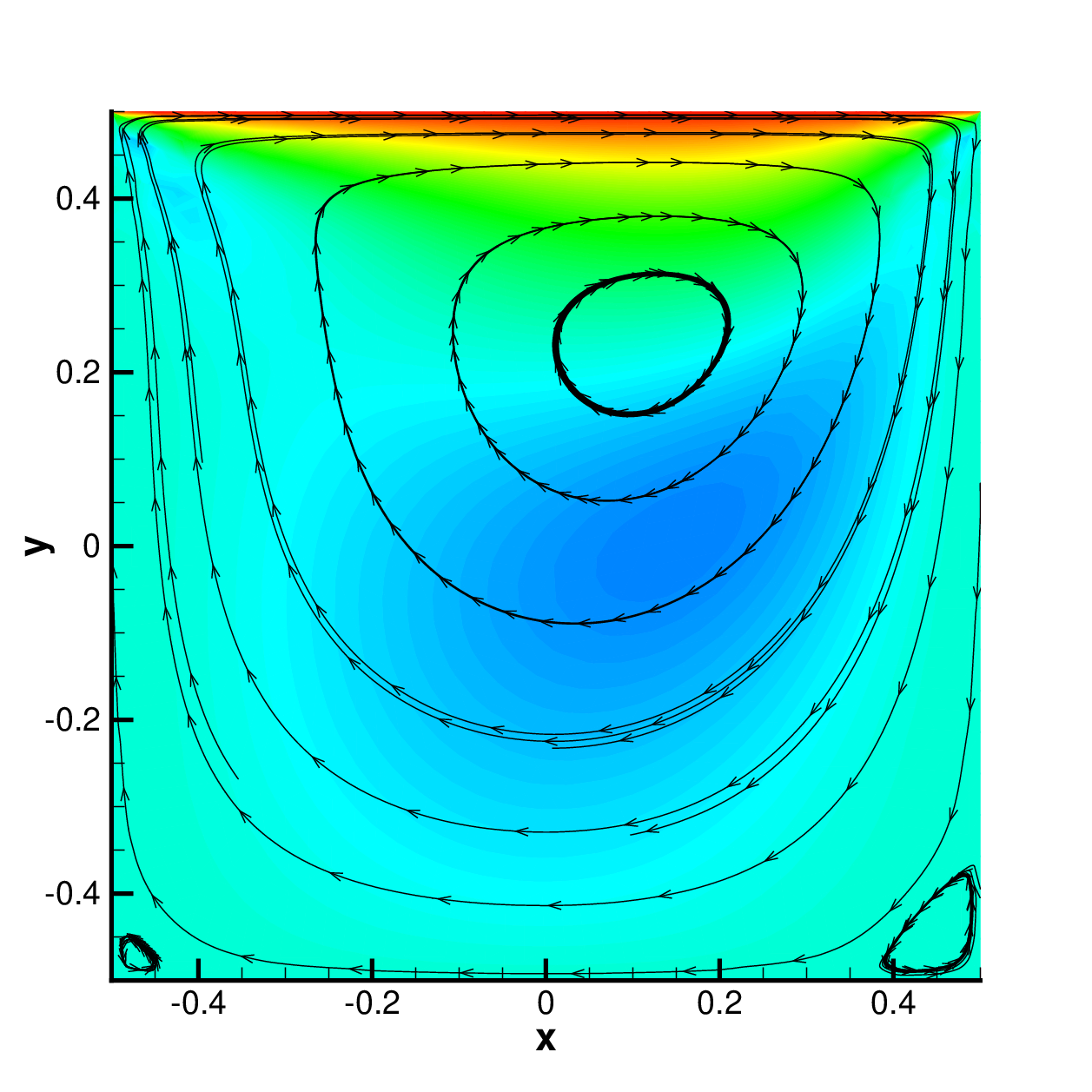} \\
    \includegraphics[width=0.48\textwidth]{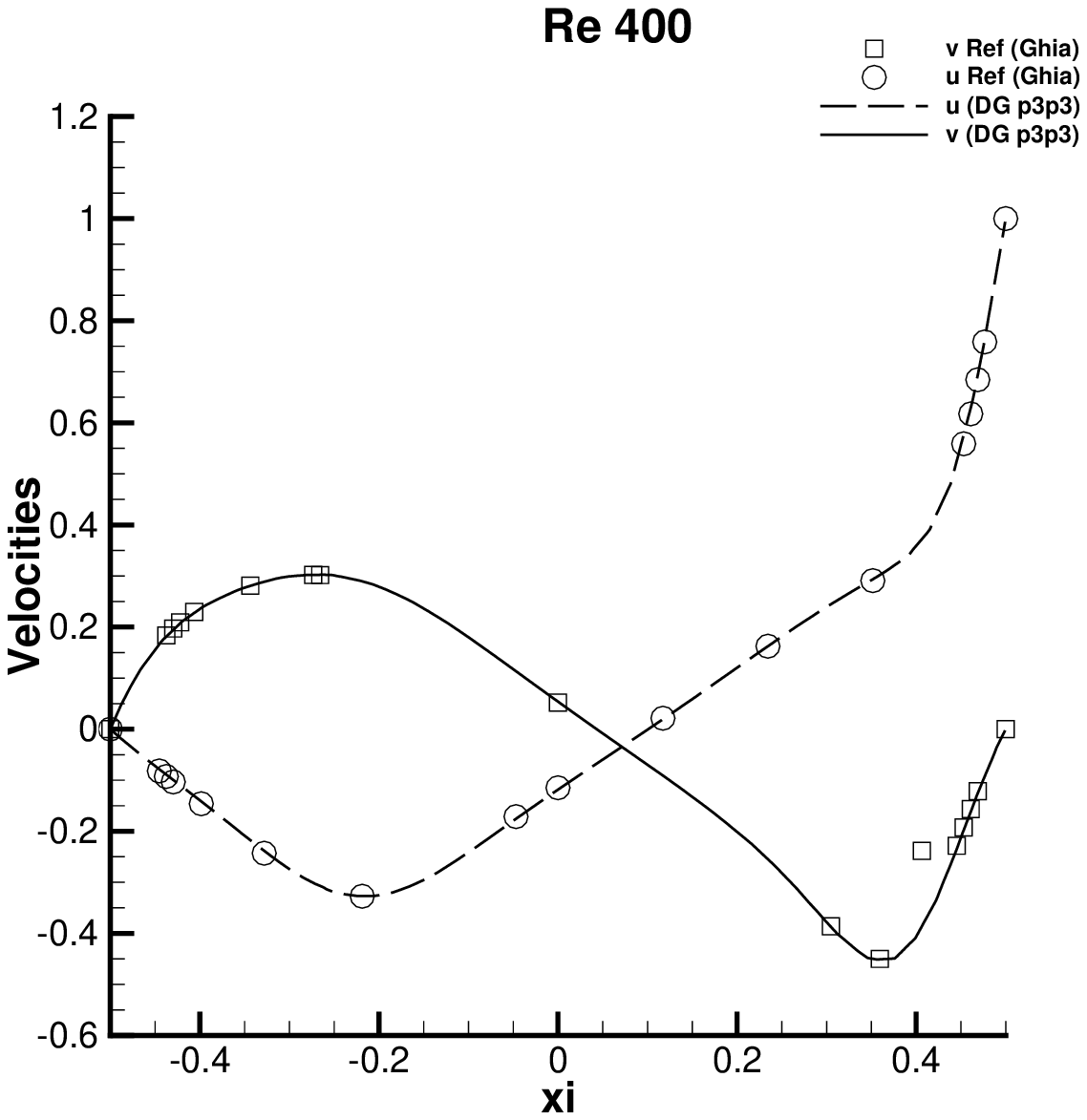}
    \includegraphics[width=0.48\textwidth]{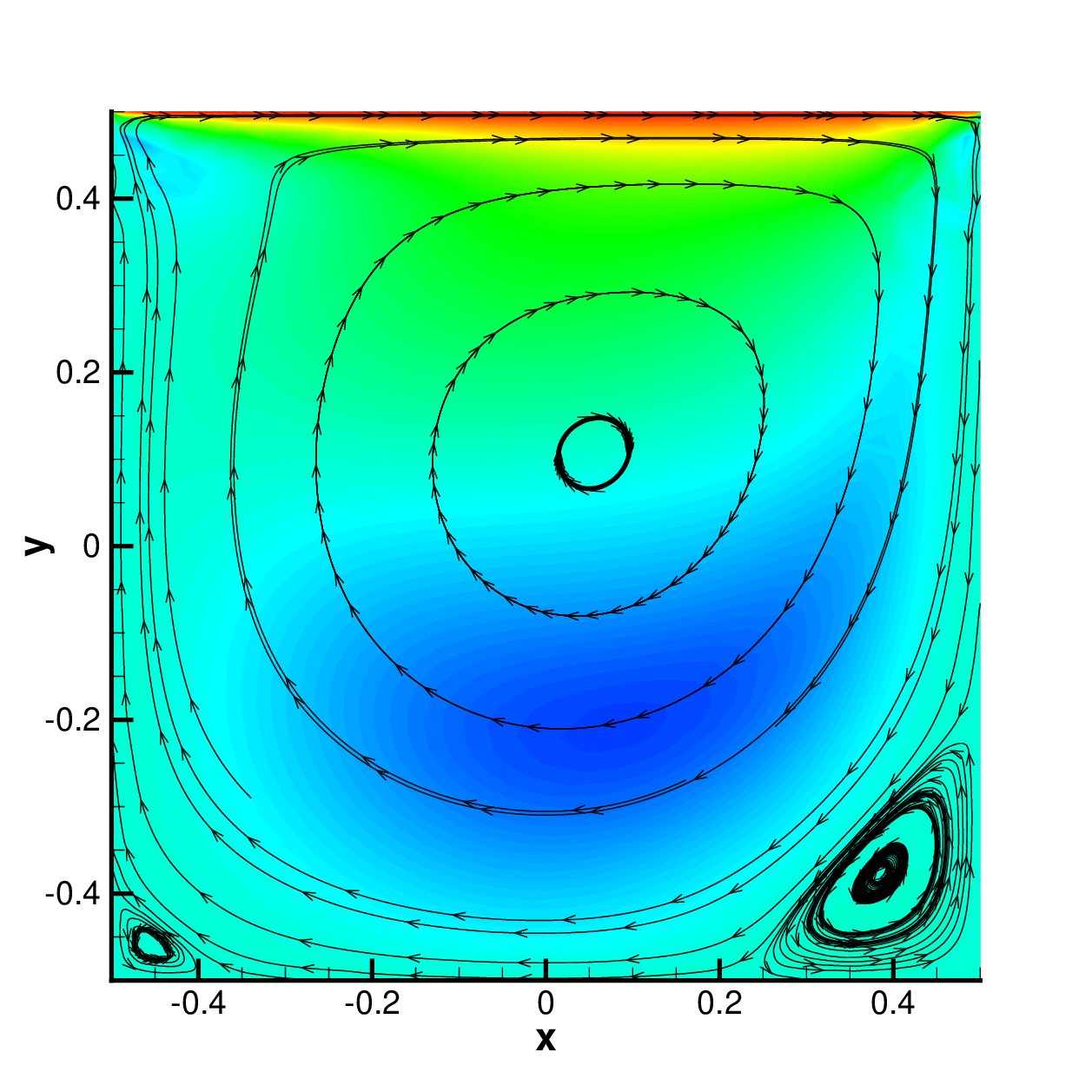} \\
    \caption{Velocity profiles (left) and streamlines (right) at Reynolds numbers Re$=100$ and Re$=400$ for the lid-driven cavity problem.}
    \label{fig.C.1}
		\end{center}
\end{figure}

\begin{figure}[ht]
    \begin{center}
    \includegraphics[width=0.48\textwidth]{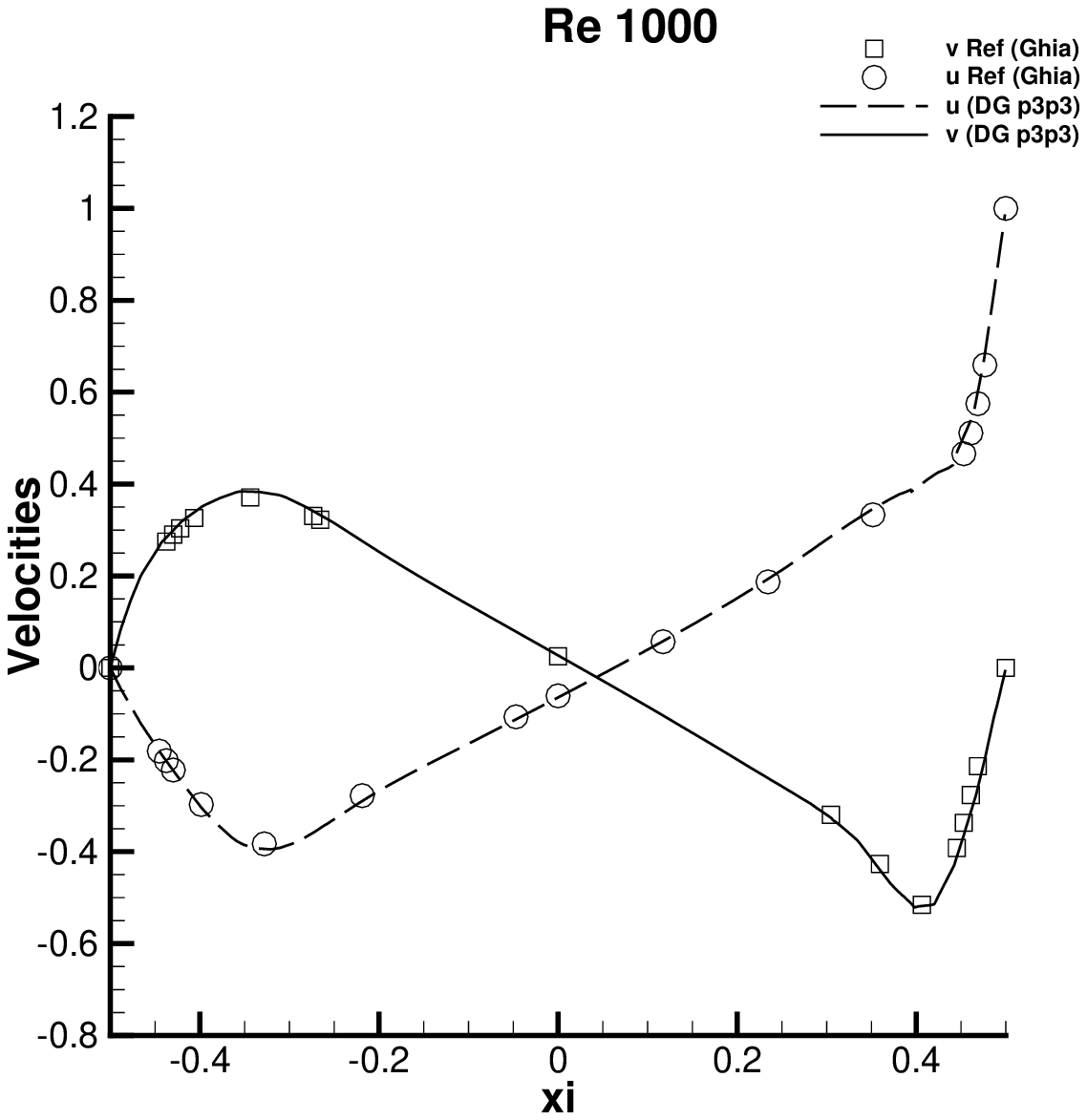}
    \includegraphics[width=0.48\textwidth]{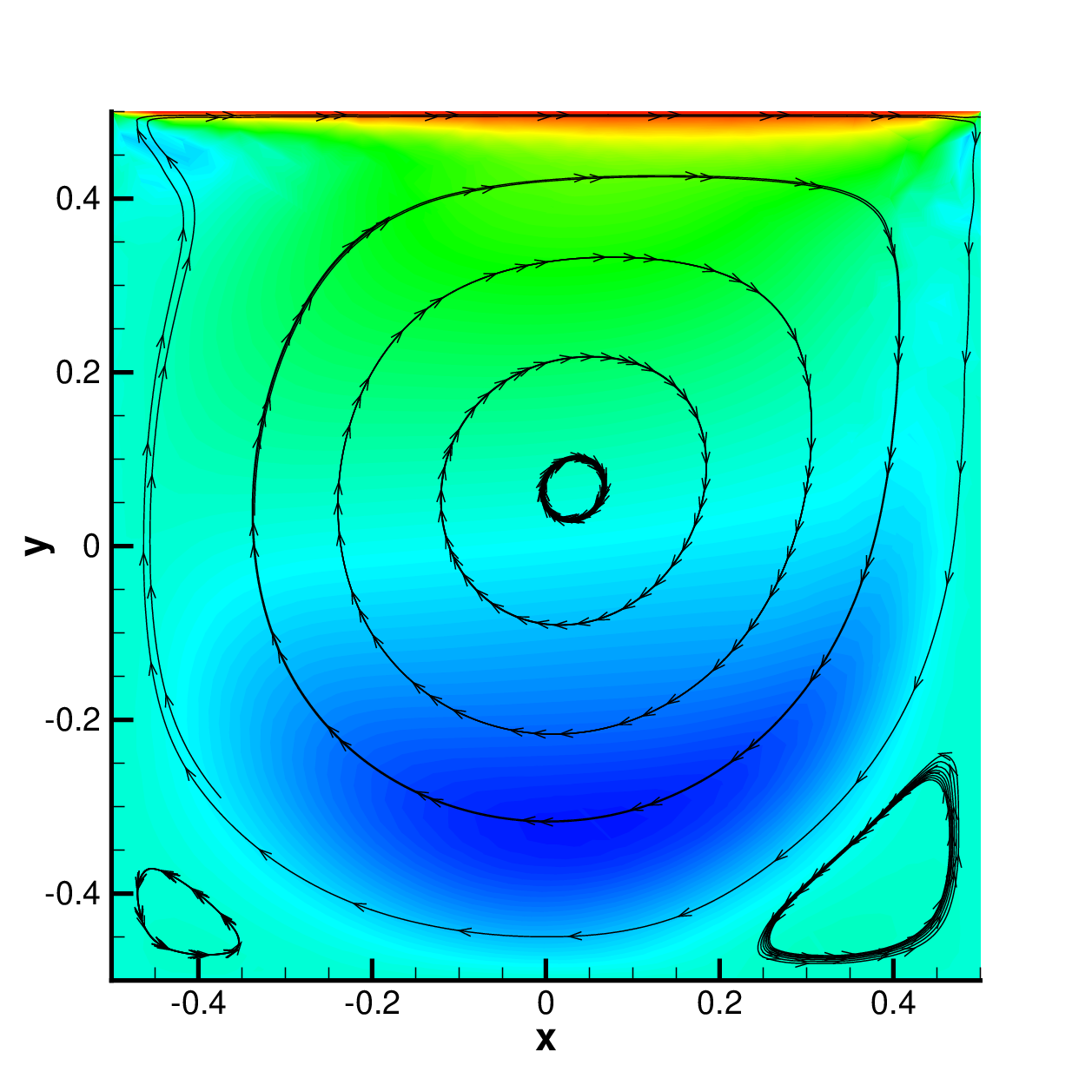} \\
    \includegraphics[width=0.48\textwidth]{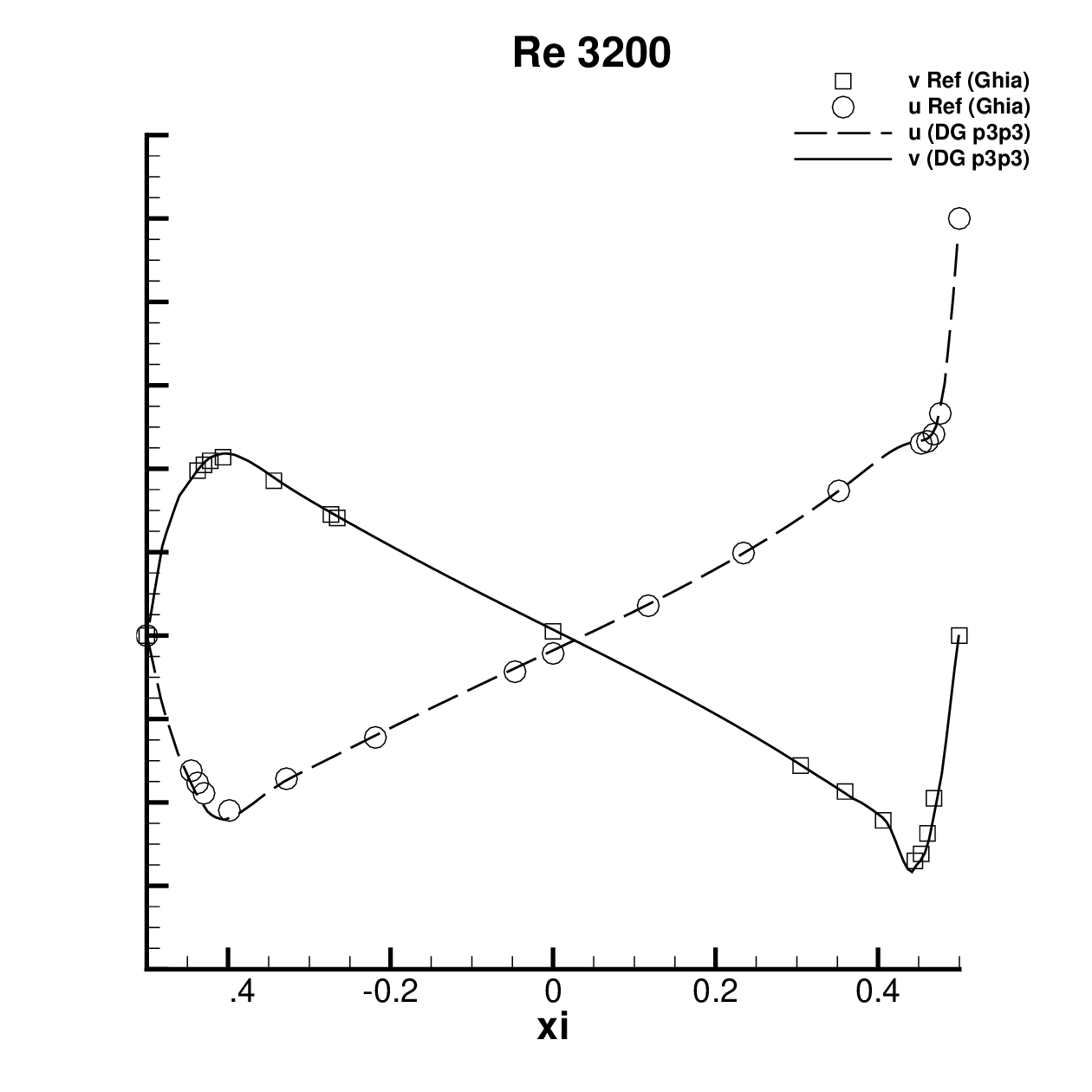}
    \includegraphics[width=0.48\textwidth]{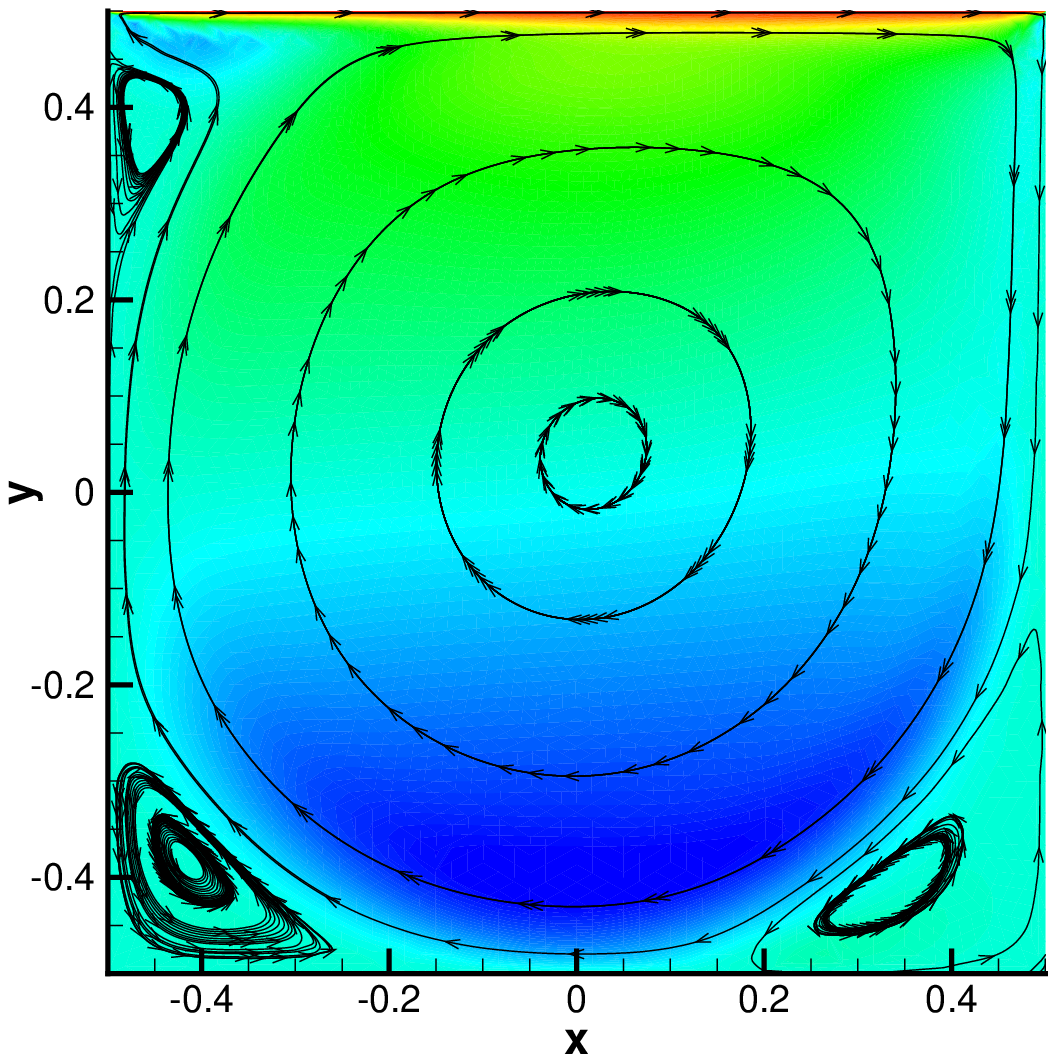} \\
    \caption{Velocity profiles (left) and streamlines (right) at Reynolds numbers Re$=1000$ and Re$=3200$ for the lid-driven cavity problem.}
    \label{fig.C.1.1}
		\end{center}
\end{figure}

For the present test $\Delta t$ is taken according to the CFL condition \eref{eq:CFLC} and $t_{end}=150$. 
According to \cite{Ghia1982}, primary and corner vortices appear from $Re=100$ to $Re=3200$, a comparison of the velocities against the data presented by 
Ghia et al in \cite{Ghia1982}, as well as the streamline plots are shown in Figures $\ref{fig.C.1}$ and $\ref{fig.C.1.1}$. A very good agreement is 
obtained in all cases, even if a very coarse grid has been used in space and time. 

\subsection{Flow over a circular cylinder}
\begin{figure}[ht]
    \begin{center}
    \includegraphics[width=0.49\textwidth]{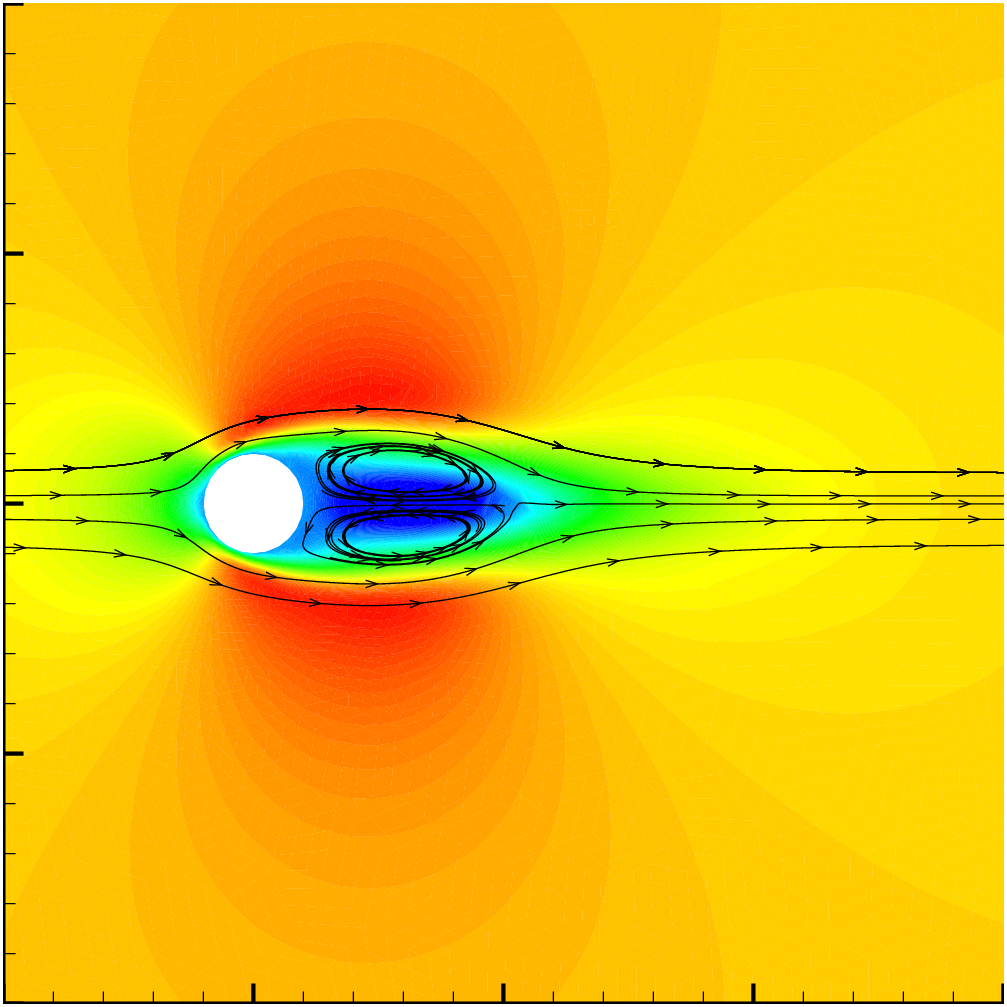}
    \includegraphics[width=0.49\textwidth]{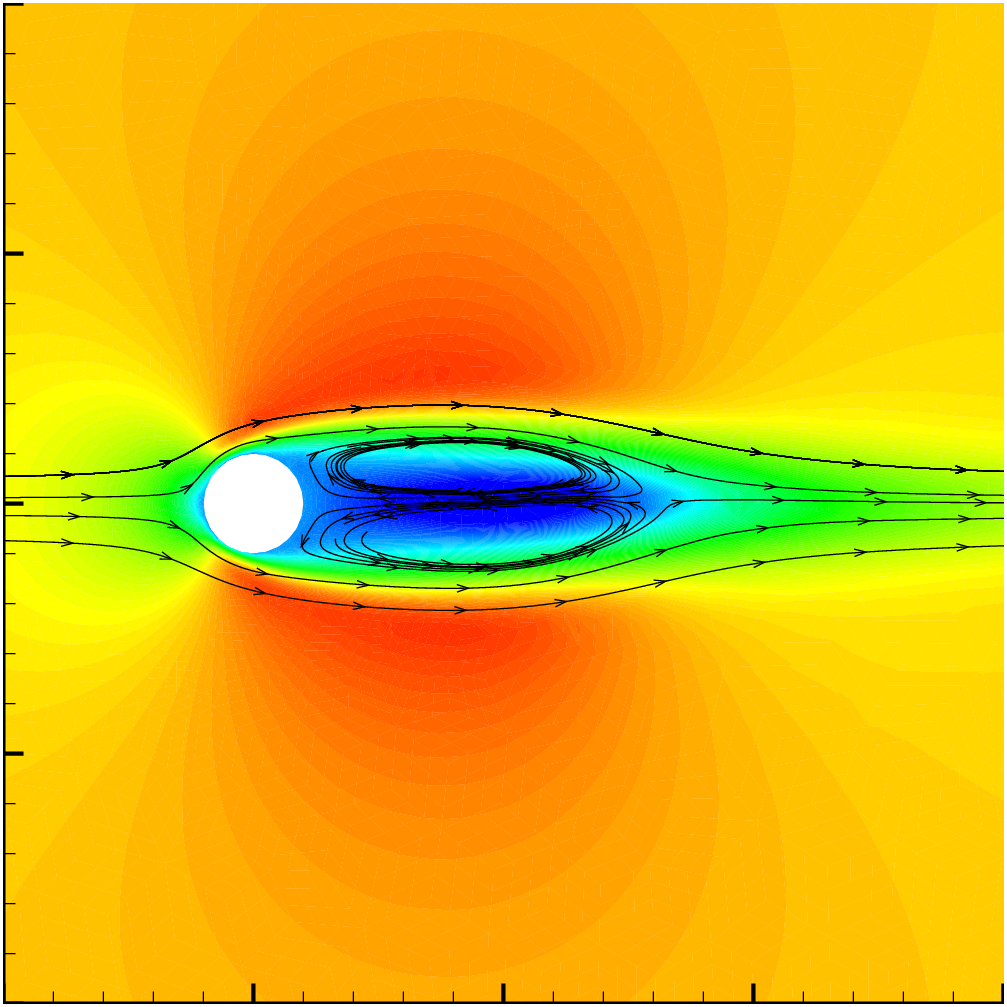} \\
    \includegraphics[width=0.49\textwidth]{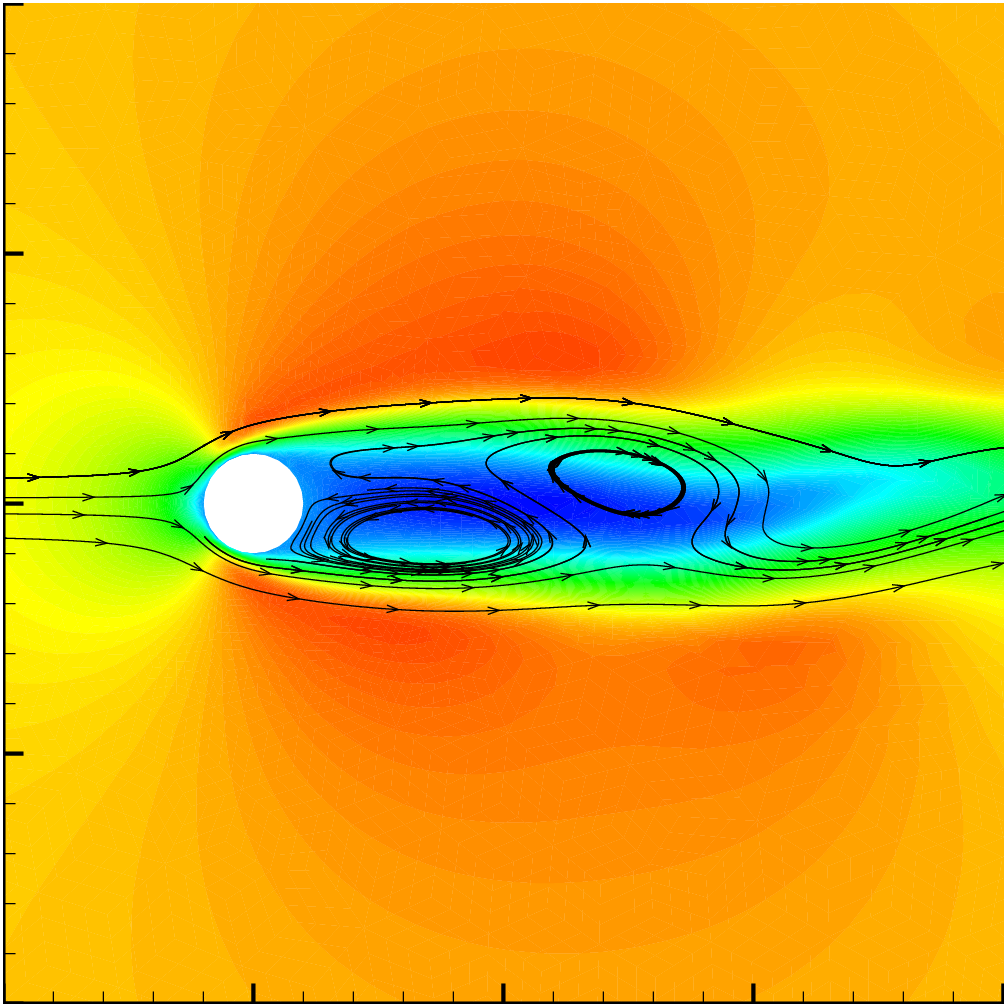}
    \includegraphics[width=0.49\textwidth]{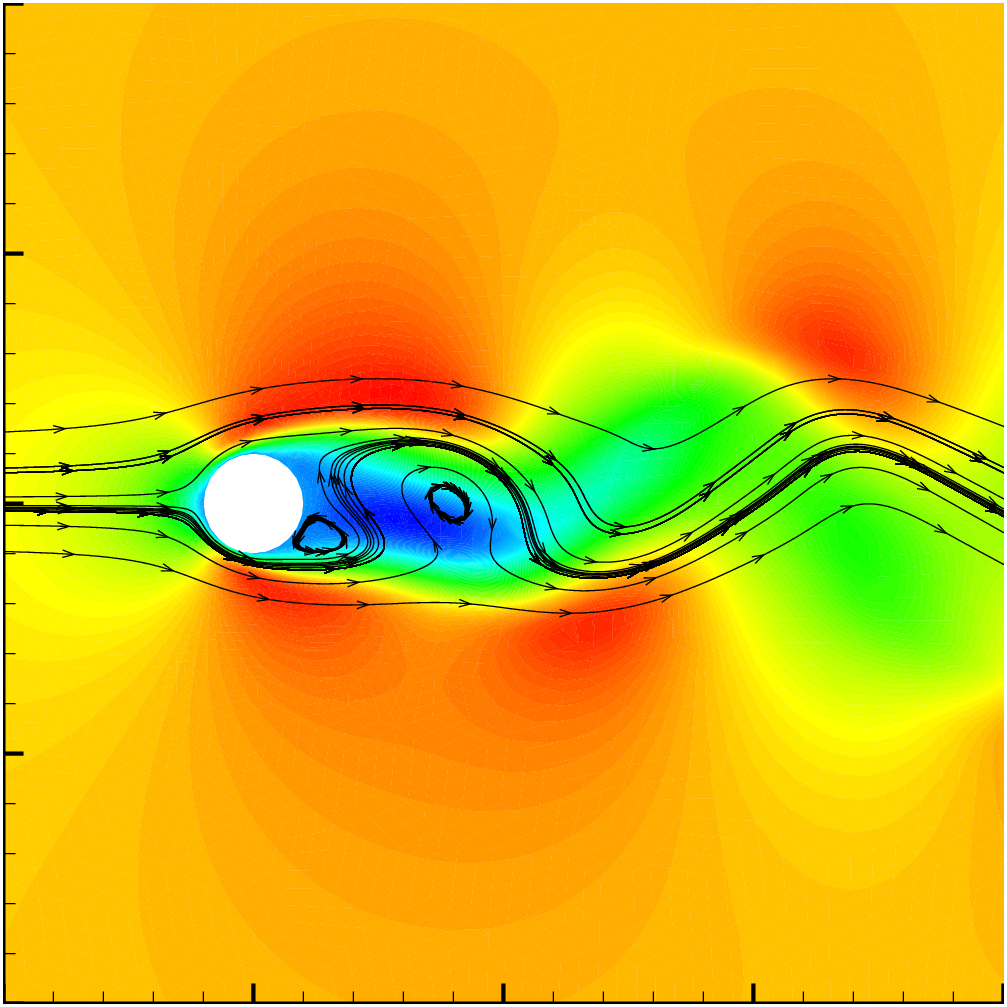}
    \caption{Streamlines along the circular cylinder at times, from top left to bottom right, $t=25,50,100$ and $t=200$}
    \label{fig.VK2}
	\end{center}
\end{figure}

\begin{figure}[ht]
    \begin{center}
    \includegraphics[width=1.0\textwidth]{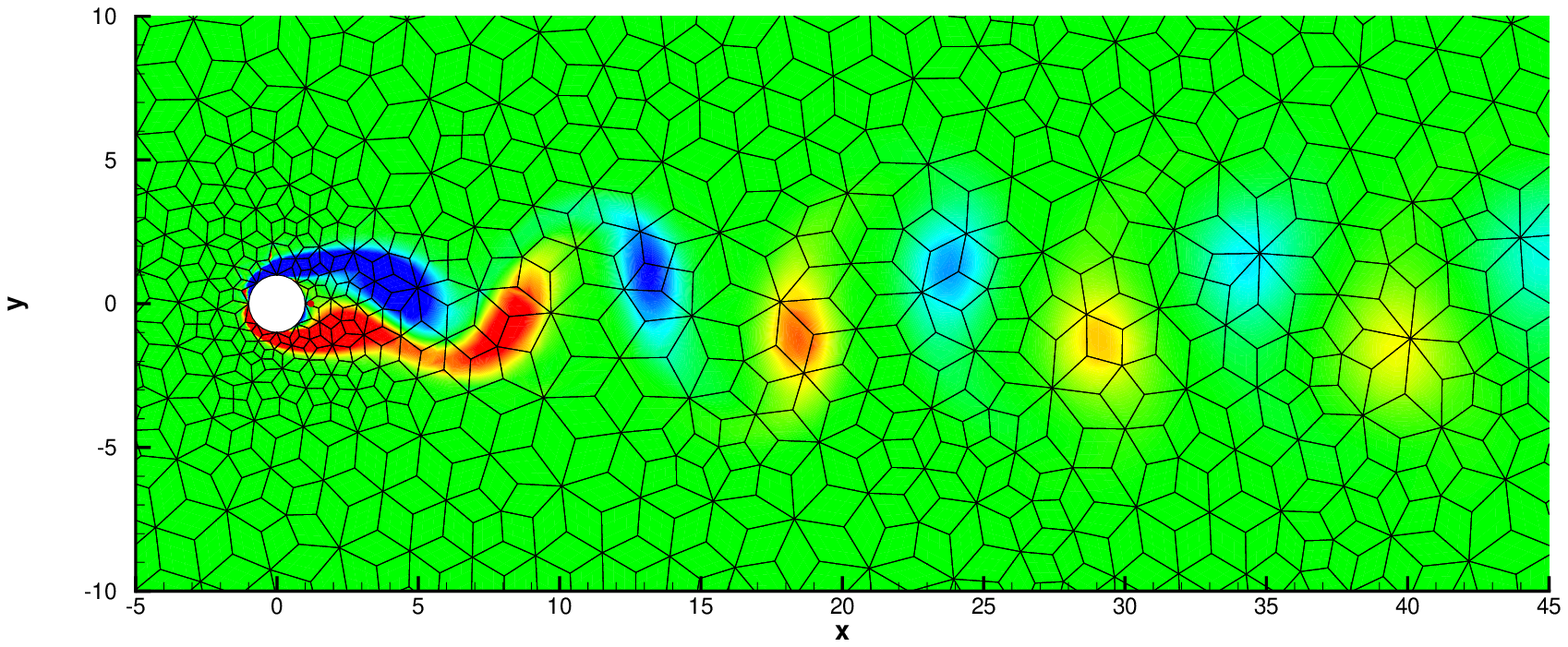} \\
    \includegraphics[width=1.0\textwidth]{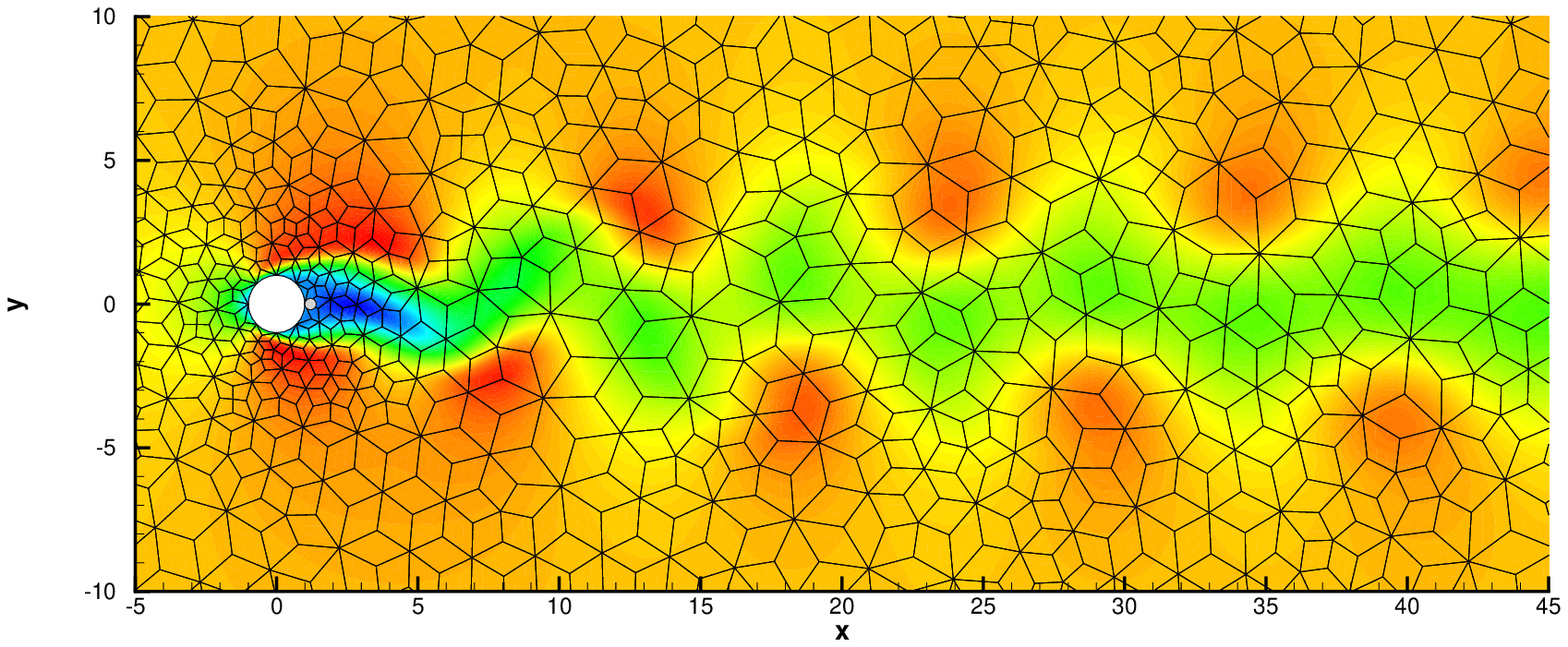}
    \caption{Laminar viscous flow past a circular cylinder. Profile for the vorticity and horizontal velocity $u$ at time $t=300$ }
    \label{fig.VK1}
		\end{center}
\end{figure}

In this section we consider the flow over a circular cylinder. In this case, the use of an isoparametric finite element approach is mandatory to represent the curved geometry of the 
cylinder wall, see \cite{BassiRebay,2DSIUSW}. 
We consider here the viscous case in order to show the formation of the von Karman vortex street. We take a sufficiently large domain $\Omega=[-20, 80]\times [-20, 20]-\{ \sqrt{x^2+y^2}\leq 1 \}$ and 
we cover it with only $N_i=1702$ triangles. Note that the chosen grid is extremely coarse compared to the dimension of the domain $\Omega$. The characteristic average size of the mesh is $h=1.295$ and 
the smallest element size is about $h_{\min}=0.347$. 
 As initial condition we set $\mathbf{v}(x,y,0)=(\bar{u},0)$, where $\bar{u}$ is the inlet velocity, taking $\bar{u}=0.5$ in our case. For the present test we use $\Delta t$ according to \eref{eq:CFLC}; 
$p=3$; $p_{\gamma}=2$. The velocity $(\bar{u},0)$ is prescribed at the left boundary while passage boundary conditions are imposed on the other external edges of the domain. Finally a no-slip wall  boundary condition is imposed on the cylinder surface. A plot of the streamlines is reported in Figure $\ref{fig.VK2}$ at several output times.
The resulting profiles for the vorticity and the horizontal velocity $u$ are plotted in Figure $\ref{fig.VK1}$, as well as the dual grid elements for $Re=100$.  
As shown in Figure $\ref{fig.VK2}$, two vortices are initially generated at the circular cylinder and then, several vortices leave the cylinder and generate the Von Karman street 
as we can see in Figure $\ref{fig.VK1}$. 
The resulting Strouhal number for the present test is $\textnormal{St}=\frac{fd}{\bar{u}}=0.1647$ that is in good agreement with $\textnormal{St}=0.1649$ obtained by Qu et al. in \cite{Lixia2013}. 

\section{Conclusions}
\label{sec.concl}
A novel high order accurate staggered semi-implicit space-time discontinuous Galerkin scheme has been proposed for the solution of the two-dimensional incompressible Navier-Stokes 
equations on unstructured curved triangular meshes. The use of a staggered grid makes our scheme different from the space-time DG schemes proposed in \cite{Rhebergen2012,Rhebergen2013}. 
The high order in space and time was verified up to $p=4$ against available exact solutions for several test cases that include a manufactured solution using source terms, 
the viscosity-dominated Womersley problem and the well-known Taylor-Green vortex problem with periodic boundary conditions.  
The numerical results agree very well with the reference data for all test cases under consideration. 
In the special case $p_\gamma=0$ the numerical method proposed in this paper reduces exactly to the semi-implicit staggered DG scheme forwarded in \cite{2STINS}, so it can be seen as 
its natural extension to high order of accuracy in time. 

Furthermore, the use of matrices that depend only on the geometry and on the polynomial degree and that hence can be precomputed before runtime, as well as a very good sparsity pattern 
involved in the solution of the main system for the pressure, leads to a computationally efficient scheme. Actually, we have solved all our test problems with a matrix-free implementation 
of the GMRES method \cite{GMRES}, without the use of any preconditioner.  

Future research will concern the extension of the scheme to the fully three-dimensional case on unstructured tetrahedral and hexahedral meshes and its application to turbulent flows. 


\section*{Acknowledgments}
The research presented in this paper was partially funded by the European Research Council (ERC) under the European Union's Seventh Framework
Programme (FP7/2007-2013) within the research project \textit{STiMulUs}, ERC Grant agreement no. 278267.
\clearpage
\bibliographystyle{elsarticle-num}
\bibliography{SIDG}

\end{document}

%% file: ugrid.tex
\begin{tikzpicture}[line cap=round,line join=round,>=triangle 45,x=0.6373937677053826cm,y=0.6177884615384613cm]
\clip(2.11,-8.53) rectangle (16.23,3.95);
\fill[color=zzttqq,fill=zzttqq,fill opacity=0.1] (5.19,-3.03) -- (9,3) -- (13,-5) -- cycle;
\fill[color=qqwuqq,fill=qqwuqq,fill opacity=0.05] (9,3) -- (14.49,2.37) -- (13,-5) -- cycle;
\fill[color=qqwuqq,fill=qqwuqq,fill opacity=0.05] (9,3) -- (4.13,2.43) -- (5.19,-3.03) -- cycle;
\fill[color=qqwuqq,fill=qqwuqq,fill opacity=0.05] (5.19,-3.03) -- (4.77,-7.83) -- (13,-5) -- cycle;
\fill[color=zzttqq,fill=zzttqq,fill opacity=0.1] (9.27,-1.79) -- (9,3) -- (5.67,0.43) -- (5.19,-3.03) -- cycle;
\fill[color=zzttqq,fill=zzttqq,fill opacity=0.1] (9.27,-1.79) -- (9,3) -- (12.15,1.37) -- (13,-5) -- cycle;
\fill[color=zzttqq,fill=zzttqq,fill opacity=0.1] (5.19,-3.03) -- (7.47,-5.49) -- (13,-5) -- (9.27,-1.79) -- cycle;
\fill[color=qqttzz,fill=qqttzz,fill opacity=0.1] (13,-5) -- (5.19,-3.03) -- (9.27,-1.79) -- cycle;
\draw [color=zzttqq] (5.19,-3.03)-- (9,3);
\draw [color=zzttqq] (9,3)-- (13,-5);
\draw [color=zzttqq] (13,-5)-- (5.19,-3.03);
\draw [color=qqwuqq] (9,3)-- (14.49,2.37);
\draw [color=qqwuqq] (14.49,2.37)-- (13,-5);
\draw [color=qqwuqq] (13,-5)-- (9,3);
\draw [color=qqwuqq] (9,3)-- (4.13,2.43);
\draw [color=qqwuqq] (4.13,2.43)-- (5.19,-3.03);
\draw [color=qqwuqq] (5.19,-3.03)-- (9,3);
\draw [color=qqwuqq] (5.19,-3.03)-- (4.77,-7.83);
\draw [color=qqwuqq] (4.77,-7.83)-- (13,-5);
\draw [color=qqwuqq] (13,-5)-- (5.19,-3.03);
\draw (9.51,-1.19) node[anchor=north west] {$i$};
\draw (12.41,1.67) node[anchor=north west] {$i_1$};
\draw (5.67,0.67) node[anchor=north west] {$i_2$};
\draw (7.47,-5.25) node[anchor=north west] {$i_3$};
\draw (11.25,-0.47) node[anchor=north west] {$j_1$};
\draw (6.95,1.43) node[anchor=north west] {$j_2$};
\draw (8.93,-3.93) node[anchor=north west] {$j_3$};
\draw (4.39,-2.79) node[anchor=north west] {$n_1$};
\draw (13.07,-4.83) node[anchor=north west] {$n_2$};
\draw (9.05,4.03) node[anchor=north west] {$n_3$};
\draw (7.71,0.35) node[anchor=north west] {$\TT_i$};
\draw [color=zzttqq] (9.27,-1.79)-- (9,3);
\draw [color=zzttqq] (9,3)-- (5.67,0.43);
\draw [color=zzttqq] (5.67,0.43)-- (5.19,-3.03);
\draw [color=zzttqq] (5.19,-3.03)-- (9.27,-1.79);
\draw [color=zzttqq] (9.27,-1.79)-- (9,3);
\draw [color=zzttqq] (9,3)-- (12.15,1.37);
\draw [color=zzttqq] (12.15,1.37)-- (13,-5);
\draw [color=zzttqq] (13,-5)-- (9.27,-1.79);
\draw [color=zzttqq] (5.19,-3.03)-- (7.47,-5.49);
\draw [color=zzttqq] (7.47,-5.49)-- (13,-5);
\draw [color=zzttqq] (13,-5)-- (9.27,-1.79);
\draw [color=zzttqq] (9.27,-1.79)-- (5.19,-3.03);
\draw (10.49,1.15) node[anchor=north west] {$\QQ_{j_1}$};
\draw [color=ffqqqq](10.07,-0.83) node[anchor=north west] {$\Gamma_{j_1}$};
\draw [line width=1.6pt,color=ffqqqq] (9,3)-- (13,-5);
\draw [color=qqttzz] (13,-5)-- (5.19,-3.03);
\draw [color=qqttzz] (5.19,-3.03)-- (9.27,-1.79);
\draw [color=qqttzz] (9.27,-1.79)-- (13,-5);
\draw [color=qqttzz](8.35,-2.17) node[anchor=north west] {$\TT_{i,j_3}$};
\draw (5.19,-3.03)-- (5.67,0.43);
\draw (5.67,0.43)-- (9,3);
\draw (9,3)-- (9.27,-1.79);
\draw (9.27,-1.79)-- (5.19,-3.03);
\draw (9,3)-- (12.15,1.37);
\draw (12.15,1.37)-- (13,-5);
\draw (13,-5)-- (9.27,-1.79);
\draw (13,-5)-- (7.47,-5.49);
\draw (7.47,-5.49)-- (5.19,-3.03);
\begin{scriptsize}
\fill [color=qqqqff] (5.19,-3.03) circle (1.5pt);
\fill [color=qqqqff] (9,3) circle (1.5pt);
\fill [color=qqqqff] (13,-5) circle (1.5pt);
\fill [color=qqqqff] (9.27,-1.79) circle (1.5pt);
\fill [color=qqqqff] (14.49,2.37) circle (1.5pt);
\fill [color=qqqqff] (4.13,2.43) circle (1.5pt);
\fill [color=qqqqff] (4.77,-7.83) circle (1.5pt);
\fill [color=qqqqff] (12.15,1.37) circle (1.5pt);
\fill [color=qqqqff] (5.67,0.43) circle (1.5pt);
\fill [color=qqqqff] (7.47,-5.49) circle (1.5pt);
\end{scriptsize}
\end{tikzpicture}

%% file: Omega_st.tex
\begin{tikzpicture}[line cap=round,line join=round,>=triangle 45,x=0.7cm,y=0.7cm]
\clip(-1.72,-8.88) rectangle (16.84,3.5);
\fill[color=zzttqq,fill=zzttqq,fill opacity=0.1] (9,-4) -- (3,-6) -- (8,-7) -- cycle;
\fill[color=ttzzqq,fill=ttzzqq,fill opacity=0.1] (9,-4) -- (7.06,-5.74) -- (8,-7) -- (11,-5) -- cycle;
\fill[color=zzttqq,fill=zzttqq,fill opacity=0.1] (9,1) -- (3,-1) -- (8,-2) -- cycle;
\fill[color=ttzzqq,fill=ttzzqq,fill opacity=0.1] (9,1) -- (7.06,-0.74) -- (8,-2) -- (11,0) -- cycle;
\fill[color=zzttqq,fill=zzttqq,fill opacity=0.1] (3,-6) -- (8,-7) -- (8,-2) -- (3,-1) -- cycle;
\fill[color=zzttqq,fill=zzttqq,fill opacity=0.1] (8,-7) -- (9,-4) -- (9,1) -- (8,-2) -- cycle;
\fill[color=ttzzqq,fill=ttzzqq,fill opacity=0.1] (8,-7) -- (11,-5) -- (11,0) -- (8,-2) -- cycle;
\fill[color=ttzzqq,fill=ttzzqq,fill opacity=0.25] (7.06,-0.74) -- (8,-2) -- (8,-7) -- (7.06,-5.74) -- cycle;
\fill[color=ttzzqq,fill=ttzzqq,fill opacity=0.2] (7.06,-5.74) -- (9,-4) -- (9,1) -- (7.06,-0.74) -- cycle;
\draw [color=zzttqq] (9,-4)-- (3,-6);
\draw [color=zzttqq] (3,-6)-- (8,-7);
\draw [color=zzttqq] (8,-7)-- (9,-4);
\draw [color=ttzzqq] (9,-4)-- (7.06,-5.74);
\draw [color=ttzzqq] (7.06,-5.74)-- (8,-7);
\draw [color=ttzzqq] (8,-7)-- (11,-5);
\draw [color=ttzzqq] (11,-5)-- (9,-4);
\draw [color=zzttqq] (9,1)-- (3,-1);
\draw [color=zzttqq] (3,-1)-- (8,-2);
\draw [color=zzttqq] (8,-2)-- (9,1);
\draw [color=ttzzqq] (9,1)-- (7.06,-0.74);
\draw [color=ttzzqq] (7.06,-0.74)-- (8,-2);
\draw [color=ttzzqq] (8,-2)-- (11,0);
\draw [color=ttzzqq] (11,0)-- (9,1);
\draw [color=zzttqq] (3,-6)-- (8,-7);
\draw [color=zzttqq] (8,-7)-- (8,-2);
\draw [color=zzttqq] (8,-2)-- (3,-1);
\draw [color=zzttqq] (3,-1)-- (3,-6);
\draw [color=zzttqq] (8,-7)-- (9,-4);
\draw [color=zzttqq] (9,-4)-- (9,1);
\draw [color=zzttqq] (9,1)-- (8,-2);
\draw [color=zzttqq] (8,-2)-- (8,-7);
\draw [color=ttzzqq] (8,-7)-- (11,-5);
\draw [color=ttzzqq] (11,-5)-- (11,0);
\draw [color=ttzzqq] (11,0)-- (8,-2);
\draw [color=ttzzqq] (8,-2)-- (8,-7);
\draw [color=ttzzqq] (7.06,-0.74)-- (8,-2);
\draw [color=ttzzqq] (8,-2)-- (8,-7);
\draw [color=ttzzqq] (8,-7)-- (7.06,-5.74);
\draw [color=ttzzqq] (7.06,-5.74)-- (7.06,-0.74);
\draw [color=ttzzqq] (7.06,-5.74)-- (9,-4);
\draw [color=ttzzqq] (9,-4)-- (9,1);
\draw [color=ttzzqq] (9,1)-- (7.06,-0.74);
\draw [color=ttzzqq] (7.06,-0.74)-- (7.06,-5.74);
\draw (5.82,-0.2) node[anchor=north west] {$\TT_i$};
\draw (9.24,0.36) node[anchor=north west] {$\QQ_j$};
\draw [->] (0,-6) -- (0.02,-4.5);
\draw [->] (0,-6) -- (1,-5);
\draw [->] (0,-6) -- (1.38,-6.26);
\draw (0.56,-4.26) node[anchor=north west] {$y$};
\draw (0.88,-6.22) node[anchor=north west] {$x$};
\draw (-0.58,-3.98) node[anchor=north west] {$t$};
\draw (1.94,-0.34) node[anchor=north west] {$t^{n+1}$};
\draw (2.22,-5.42) node[anchor=north west] {$t^{n}$};
\draw (5.06,-2.82) node[anchor=north west] {$\TT_i^{st}$};
\draw (9.44,-2.44) node[anchor=north west] {$\QQ_j^{st}$};
\draw (9,1)-- (15,0);
\draw (15,0)-- (8,-2);
\draw (15,-5)-- (15,0);
\draw (15,-5)-- (9,-4);
\draw (15,-5)-- (8,-7);
\end{tikzpicture}

%% file: Layer.tex
\begin{tikzpicture}[line cap=round,line join=round,>=triangle 45,x=0.7cm,y=0.7cm]
\clip(1.91,-6.69) rectangle (11.95,3.72);
\fill[color=zzttqq,fill=zzttqq,fill opacity=0.1] (9,-2) -- (3,-4) -- (8,-5) -- cycle;
\fill[color=qqttcc,fill=qqttcc,fill opacity=0.1] (9,-1) -- (3,-3) -- (8,-4) -- cycle;
\fill[color=qqwwtt,fill=qqwwtt,fill opacity=0.1] (9,1) -- (3,-1) -- (8,-2) -- cycle;
\fill[color=zzttqq,fill=zzttqq,fill opacity=0.1] (9,2) -- (3,0) -- (8,-1) -- cycle;
\fill[color=zzttqq,fill=zzttqq,fill opacity=0.1] (3,-4) -- (3,0) -- (8,-1) -- (8,-5) -- cycle;
\fill[color=zzttqq,fill=zzttqq,fill opacity=0.1] (8,-5) -- (9,-2) -- (9,2) -- (8,-1) -- cycle;
\draw [color=zzttqq] (9,-2)-- (3,-4);
\draw [color=zzttqq] (3,-4)-- (8,-5);
\draw [color=zzttqq] (8,-5)-- (9,-2);
\draw [color=qqttcc] (9,-1)-- (3,-3);
\draw [color=qqttcc] (3,-3)-- (8,-4);
\draw [color=qqttcc] (8,-4)-- (9,-1);
\draw [color=qqwwtt] (9,1)-- (3,-1);
\draw [color=qqwwtt] (3,-1)-- (8,-2);
\draw [color=qqwwtt] (8,-2)-- (9,1);
\draw [color=zzttqq] (9,2)-- (3,0);
\draw [color=zzttqq] (3,0)-- (8,-1);
\draw [color=zzttqq] (8,-1)-- (9,2);
\draw [color=zzttqq] (3,-4)-- (3,0);
\draw [color=zzttqq] (3,0)-- (8,-1);
\draw [color=zzttqq] (8,-1)-- (8,-5);
\draw [color=zzttqq] (8,-5)-- (3,-4);
\draw [color=zzttqq] (8,-5)-- (9,-2);
\draw [color=zzttqq] (9,-2)-- (9,2);
\draw [color=zzttqq] (9,2)-- (8,-1);
\draw [color=zzttqq] (8,-1)-- (8,-5);
\draw [color=qqttcc](9.8,-1.8) node[anchor=north west] {$\mathbf{\ell_2=1}$};
\draw [color=qqwwtt](9.8,0.43) node[anchor=north west] {\parbox{1.39 cm}{$\mathbf{\ell_2 \\ =2}$}};
\draw [->,color=qqttcc] (9.55,-2.38) -- (8.53,-2.4);
\draw [->,color=qqwwtt] (9.55,-0.19) -- (8.6,-0.2);
\draw [->] (9,3) -- (9,2);
\draw [->] (8,-6) -- (8,-5);
\draw [->] (3,1) -- (3,0);
\draw (2.26,1.79) node[anchor=north west] {$\ell_1=1$};
\draw (7.26,-5.5) node[anchor=north west] {$\ell_1=2$};
\draw (8.29,3.82) node[anchor=north west] {$\ell_1=3$};
\begin{scriptsize}
\draw [color=qqqqff] (9,-1)-- ++(-2.5pt,-2.5pt) -- ++(5.0pt,5.0pt) ++(-5.0pt,0) -- ++(5.0pt,-5.0pt);
\draw [color=qqqqff] (3,-3)-- ++(-2.5pt,-2.5pt) -- ++(5.0pt,5.0pt) ++(-5.0pt,0) -- ++(5.0pt,-5.0pt);
\draw [color=qqqqff] (8,-4)-- ++(-2.5pt,-2.5pt) -- ++(5.0pt,5.0pt) ++(-5.0pt,0) -- ++(5.0pt,-5.0pt);
\draw [color=qqwwtt] (9,1)-- ++(-2.5pt,-2.5pt) -- ++(5.0pt,5.0pt) ++(-5.0pt,0) -- ++(5.0pt,-5.0pt);
\draw [color=qqwwtt] (3,-1)-- ++(-2.5pt,-2.5pt) -- ++(5.0pt,5.0pt) ++(-5.0pt,0) -- ++(5.0pt,-5.0pt);
\draw [color=qqwwtt] (8,-2)-- ++(-2.5pt,-2.5pt) -- ++(5.0pt,5.0pt) ++(-5.0pt,0) -- ++(5.0pt,-5.0pt);
\end{scriptsize}
\end{tikzpicture}

%% file: SpaceTimeINS.bbl
\begin{thebibliography}{10}
\expandafter\ifx\csname url\endcsname\relax
  \def\url#1{\texttt{#1}}\fi
\expandafter\ifx\csname urlprefix\endcsname\relax\def\urlprefix{URL }\fi
\expandafter\ifx\csname href\endcsname\relax
  \def\href#1#2{#2} \def\path#1{#1}\fi

\bibitem{markerandcell}
F.~Harlow, J.~Welch, Numerical calculation of time-dependent viscous
  incompressible flow of fluid with a free surface, Physics of Fluids 8 (1965)
  2182--2189.

\bibitem{patankarspalding}
V.~Patankar, B.~Spalding, A calculation procedure for heat, mass and momentum
  transfer in three-dimensional parabolic flows, International Journal of Heat
  and Mass Transfer 15 (1972) 1787--1806.

\bibitem{patankar}
V.~Patankar, Numerical {Heat} {Transfer} and {Fluid} {Flow}, Hemisphere
  Publishing Corporation, 1980.

\bibitem{vanKan}
J.~van Kan, {A second-order accurate pressure correction method for viscous
  incompressible flow}, SIAM Journal on Scientific and Statistical Computing 7
  (1986) 870--891.

\bibitem{TaylorHood}
C.~Taylor, P.~Hood, {A numerical solution of the Navier-Stokes equations using
  the finite element technique}, Computers and Fluids 1 (1973) 73--100.

\bibitem{SUPG}
A.~Brooks, T.~Hughes, {Stream-line upwind/Petrov Galerkin formulstion for
  convection dominated flows with particular emphasis on the incompressible
  Navier-Stokes equation}, Computer Methods in Applied Mechanics and
  Engineering 32 (1982) 199--259.

\bibitem{SUPG2}
T.~Hughes, M.~Mallet, M.~Mizukami, {A new finite element formulation for
  computational fluid dynamics: II. Beyond SUPG}, Computer Methods in Applied
  Mechanics and Engineering 54 (1986) 341--–355.

\bibitem{Fortin}
M.~Fortin, {Old and new finite elements for incompressible flows},
  International Journal for Numerical Methods in Fluids 1 (1981) 347--364.

\bibitem{Verfuerth}
R.~Verf\"urth, {Finite element approximation of incompressible Navier-Stokes
  equations with slip boundary condition II}, Numerische Mathematik 59 (1991)
  615--636.

\bibitem{Rannacher1}
J.~G. Heywood, R.~Rannacher, {Finite element approximation of the nonstationary
  Navier-Stokes Problem. I. Regularity of solutions and second order error
  estimates for spatial discretization}, SIAM Journal on Numerical Analysis 19
  (1982) 275--311.

\bibitem{Rannacher3}
J.~G. Heywood, R.~Rannacher, {Finite element approximation of the nonstationary
  Navier-Stokes Problem. III. Smoothing property and higher order error
  estimates for spatial discretization}, SIAM Journal on Numerical Analysis 25
  (1988) 489--512.

\bibitem{Bassi2007}
F.~Bassi, A.~Crivellini, D.~D. Pietro, S.~Rebay, {An implicit high-order
  discontinuous Galerkin method for steady and unsteady incompressible flows},
  Computers and Fluids 36 (2007) 1529--1546.

\bibitem{Shahbazi2007}
K.~Shahbazi, P.~F. Fischer, C.~R. Ethier, A high-order discontinuous galerkin
  method for the unsteady incompressible navier-stokes equations, Journal of
  Computational Physics 222 (2007) 391--407.

\bibitem{Ferrer2011}
E.~Ferrer, R.~Willden, A high order discontinuous galerkin finite element
  solver for the incompressible navier–stokes equations, Computer and Fluids 46
  (2011) 224--230.

\bibitem{Nguyen2011}
N.~Nguyen, J.~Peraire, B.~Cockburn, An implicit high-order hybridizable
  discontinuous galerkin method for the incompressible navier-stokes equations,
  Journal of Computational Physics 230 (2011) 1147--1170.

\bibitem{Rhebergen2012}
S.~Rhebergen, B.~Cockburn, {A space–time hybridizable discontinuous Galerkin
  method for incompressible flows on deforming domains}, Journal of
  Computational Physics 231 (2012) 4185--4204.

\bibitem{Rhebergen2013}
S.~Rhebergen, B.~Cockburn, J.~J. van~der Vegt, {A space–time discontinuous
  Galerkin method for the incompressible Navier–Stokes equations}, Journal of
  Computational Physics 233 (2013) 339--358.

\bibitem{Crivellini2013}
A.~Crivellini, V.~D'Alessandro, F.~Bassi, {High-order discontinuous Galerkin
  solutions of three-dimensional incompressible RANS equations}, Computers and
  Fluids 81 (2013) 122--133.

\bibitem{KleinKummerOberlack2013}
B.~Klein, F.~Kummer, M.~Oberlack, {A SIMPLE based discontinuous Galerkin solver
  for steady incompressible flows}, Journal of Computational Physics 237 (2013)
  235--250.

\bibitem{Bassi2006}
F.~Bassi, A.~Crivellini, D.~D. Pietro, S.~Rebay, On a robust discontinuous
  galerkin technique for the solution of compressible flow, Journal of
  Computational Physics 218 (2006) 208--221.

\bibitem{chorin1}
A.~Chorin, A numerical method for solving incompressible viscous flow problems,
  Journal of Computational Physics 2 (1967) 12--26.

\bibitem{chorin2}
A.~Chorin, Numerical solution of the {Navier--Stokes} equations, Mathematics of
  Computation 23 (1968) 341--354.

\bibitem{Baumann1999311}
C.~Baumann, J.~Oden, A discontinuous hp finite element method for
  convection-diffusion problems, Computer Methods in Applied Mechanics and
  Engineering 175~(3-4) (1999) 311--341.

\bibitem{Baumann199979}
C.~Baumann, J.~Oden, A discontinuous hp finite element method for the euler and
  navier-stokes equations, International Journal for Numerical Methods in
  Fluids 31~(1) (1999) 79--95.

\bibitem{ArnoldBrezzi}
D.~Arnold, F.~Brezzi, B.~Cockburn, L.~Marini, Unified analysis of discontinuous
  galerkin methods for elliptic problems, SIAM Journal on Numerical Analysis 39
  (2002) 1749--1779.

\bibitem{TumoloBonaventuraRestelli}
G.~Tumolo, L.~Bonaventura, M.~Restelli, {A semi-implicit, semi-Lagrangian,
  p-adaptive discontinuous Galerkin method for the shallow water equations },
  Journal of Computational Physics 232 (2013) 46--67.

\bibitem{GiraldoRestelli}
F.~X. Giraldo, M.~Restelli, High-order semi-implicit time-integrators for a
  triangular discontinuous galerkin oceanic shallow water model, International
  Journal for Numerical Methods in Fluids 63 (2010) 1077--1102.

\bibitem{Dolejsi1}
V.~Dolejsi, Semi-implicit interior penalty discontinuous galerkin methods for
  viscous compressible flows, Communications in Computational Physics 4 (2008)
  231--274.

\bibitem{Dolejsi2}
V.~Dolejsi, M.~Feistauer, A semi-implicit discontinuous galerkin finite element
  method for the numerical solution of inviscid compressible flow, Journal of
  Computational Physics 198 (2004) 727--746.

\bibitem{Dolejsi3}
V.~Dolejsi, M.~Feistauer, J.~Hozman, Analysis of semi-implicit dgfem for
  nonlinear convection-diffusion problems on nonconforming meshes, Computer
  Methods in Applied Mechanics and Engineering 196 (2007) 2813--2827.

\bibitem{CentralDG1}
Y.~J. Liu, C.~W. Shu, E.~Tadmor, M.~Zhang, Central discontinuous galerkin
  methods on overlapping cells with a non-oscillatory hierarchical
  reconstruction, SIAM Journal on Numerical Analysis 45 (2007) 2442--2467.

\bibitem{CentralDG2}
Y.~J. Liu, C.~W. Shu, E.~Tadmor, M.~Zhang, L2-stability analysis of the central
  discontinuous galerkin method and a comparison between the central and
  regular discontinuous galerkin methods, Mathematical Modeling and Numerical
  Analysis 42 (2008) 593--607.

\bibitem{StaggeredDGCE1}
E.~Chung, B.~Engquist, {Optimal discontinuous Galerkin methods for wave
  propagation}, SIAM Journal on Numerical Analysis 44 (2006) 2131--2158.

\bibitem{StaggeredDG}
E.~T. Chung, C.~S. Lee, {A staggered discontinuous Galerkin method for the
  convection--diffusion equation}, Journal of Numerical Mathematics 20 (2012)
  1--31.

\bibitem{DumbserCasulli}
M.~Dumbser, V.~Casulli, A staggered semi-implicit spectral discontinuous
  galerkin scheme for the shallow water equations, Applied Mathematics and
  Computation 219~(15) (2013) 8057--8077.

\bibitem{2DSIUSW}
M.~Tavelli, M.~Dumbser, A high order semi-implicit discontinuous galerkin
  method for the two dimensional shallow water equations on staggered
  unstructured meshes, Applied Mathematics and Computation 234 (2014) 623--644.

\bibitem{2STINS}
M.~Tavelli, M.~Dumbser, A staggered arbitrary high order semi-implicit
  discontinuous galerkin method for the two dimensional incompressible
  navier-stokes equations, Applied Mathematics and Computation 248 (2014)
  70--92.

\bibitem{HirtNichols}
C.~W. Hirt, B.~D. Nichols, Volume of fluid ({VOF}) method for dynamics of free
  boundaries, Journal of Computational Physics 39 (1981) 201--225.

\bibitem{CasulliCheng1992}
V.~Casulli, R.~T. Cheng, Semi-implicit finite difference methods for
  three--dimensional shallow water flow, International Journal for Numerical
  Methods in Fluids 15 (1992) 629--648.

\bibitem{Casulli1999}
V.~Casulli, A semi-implicit finite difference method for non-hydrostatic
  free-surface flows, International Journal for Numerical Methods in Fluids 30
  (1999) 425--440.

\bibitem{CasulliWalters2000}
V.~Casulli, R.~A. Walters, An unstructured grid, three--dimensional model based
  on the shallow water equations, International Journal for Numerical Methods
  in Fluids 32 (2000) 331--348.

\bibitem{Casulli2009}
V.~Casulli, A high-resolution wetting and drying algorithm for free-surface
  hydrodynamics, International Journal for Numerical Methods in Fluids 60
  (2009) 391--408.

\bibitem{CasulliStelling2011}
V.~Casulli, G.~S. Stelling, Semi-implicit subgrid modelling of
  three-dimensional free-surface flows, International Journal for Numerical
  Methods in Fluids 67 (2011) 441--449.

\bibitem{CasulliVOF}
V.~Casulli, {A semi--implicit numerical method for the free--surface
  Navier--Stokes equations}, International Journal for Numerical Methods in
  Fluids 74 (2014) 605--622.

\bibitem{klein}
R.~Klein, N.~Botta, T.~Schneider, C.~Munz, S.Roller, A.~Meister, L.~Hoffmann,
  T.~Sonar, Asymptotic adaptive methods for multi-scale problems in fluid
  mechanics, Journal of Engineering Mathematics 39 (2001) 261--343.

\bibitem{RoMu99}
S.~Roller, C.~Munz, A low mach number scheme based on multi-scale asymptotics,
  Computing and Visualization in Science 3 (2000) 85--91.

\bibitem{ParkMPV}
J.~Park, C.-D. Munz, Multiple pressure variables methods for fluid flow at all
  mach numbers, International Journal for Numerical Methods in Fluids 49~(8)
  (2005) 905--931.

\bibitem{cbs4}
B.~Cockburn, C.~W. Shu, The {Runge}-{Kutta} discontinuous {Galerkin} method for
  conservation laws {V}: multidimensional systems, Journal of Computational
  Physics 141 (1998) 199--224.

\bibitem{CBS-convection-diffusion}
B.~Cockburn, C.~W. Shu, The local discontinuous {Galerkin} method for
  time-dependent convection diffusion systems, SIAM Journal on Numerical
  Analysis 35 (1998) 2440--2463.

\bibitem{CBS-convection-dominated}
B.~Cockburn, C.~W. Shu, {Runge}-{Kutta} discontinuous {Galerkin} methods for
  convection-dominated problems, Journal of Scientific Computing 16 (2001)
  173--261.

\bibitem{Rusanov:1961a}
V.~V. Rusanov, {Calculation of Interaction of Non--Steady Shock Waves with
  Obstacles}, J. Comput. Math. Phys. USSR 1 (1961) 267--279.

\bibitem{MunzDiffusionFlux}
G.~Gassner, F.~L\"orcher, C.~D. Munz, A contribution to the construction of
  diffusion fluxes for finite volume and discontinuous {Galerkin} schemes,
  Journal of Computational Physics 224 (2007) 1049--1063.

\bibitem{Bermudez1998}
A.~Bermudez, A.~Dervieux, J.~Desideri, M.~Vazquez, Upwind schemes for the
  two--dimensional shallow water equations with variable depth using
  unstructured meshes, Computer Methods in Applied Mechanics and Engineering
  155 (1998) 49--72.

\bibitem{Bermudez2014}
A.~Berm\'udez, J.~Ferr\'in, L.~Saavedra, M.~V\'azquez-Cend\'on, {A projection
  hybrid finite volume/element method for low-Mach number flows}, Journal of
  Computational Physics 271 (2014) 360--378.

\bibitem{USFORCE}
E.~F. Toro, A.~Hidalgo, M.~Dumbser, {FORCE} schemes on unstructured meshes {I}:
  Conservative hyperbolic systems, Journal of Computational Physics 228 (2009)
  3368--–3389.

\bibitem{ADERNSE}
M.~Dumbser, Arbitrary high order {PNPM} schemes on unstructured meshes for the
  compressible {Navier--Stokes} equations, Computers \& Fluids 39 (2010)
  60--76.

\bibitem{Dumbser2008}
M.~Dumbser, D.~S. Balsara, E.~F. Toro, C.~D. Munz, A unified framework for the
  construction of one-step finite-volume and discontinuous {Galerkin} schemes,
  Journal of Computational Physics 227 (2008) 8209--–8253.

\bibitem{CasulliZanolli}
V.~Casulli, P.~Zanolli, High resolution methods for multidimensional
  advection--diffusion problems in free--surface hydrodynamics, Ocean Modelling
  10 (2005) 137--151.

\bibitem{TavelliDumbserCasulli}
M.~Tavelli, M.~Dumbser, V.~Casulli, High resolution methods for scalar
  transport problems in compliant systems of arteries, Applied Numerical
  Mathematics 74 (2013) 62--82.

\bibitem{GMRES}
Y.~Saad, M.~Schultz, {GMRES:} a generalized minimal residual algorithm for
  solving nonsymmetric linear systems, SIAM Journal on Scientific and
  Statistical Computing 7 (1986) 856–--869.

\bibitem{Womersley1995}
J.~Womersley, Method for the calculation of velocity, rate of flow and viscous
  drag in arteries when the pressure gradient is known, Journal of Physiology
  127 (1955) 553--563.

\bibitem{Bell1989}
J.~B. Bell, P.~Coletta, H.~M. Glaz, A second-order projection method for the
  incompressible navier-stokes equations, Journal of Computational Physics 85
  (1989) 257--283.

\bibitem{Ghia1982}
U.~Ghia, K.~N. Ghia, C.~T. Shin, High-re solutions for incompressible flow
  using navier-stokes equations and multigrid method, Journal of Computational
  Physics 48 (1982) 387--411.

\bibitem{BassiRebay}
F.~Bassi, S.~Rebay, A high-order accurate discontinuous finite element method
  for the numerical solution of the compressible {Navier}-{Stokes} equations,
  Journal of Computational Physics 131 (1997) 267--279.

\bibitem{Lixia2013}
L.~Qu, C.~Norberg, L.~Davidson, S.~Peng, F.~Wang, Quantitative numerical
  analysis of flow past a circular cylinder at reynolds number between 50 and
  200, Journal of Fluids and Structures 39 (2013) 347--370.

\end{thebibliography}
